\documentclass[a4paper,10pt]{article}

\usepackage{amsthm}

\usepackage{amsmath}

\usepackage{amsfonts}

\usepackage{amssymb}

\usepackage{amsmath,amsthm,amscd,amssymb}

\usepackage{latexsym}

\usepackage{graphicx}

\usepackage[all]{xy}

\usepackage{setspace}

\theoremstyle{plain}

\newtheorem{thm}{\textbf{Theorem}}

\newtheorem{cn}{\textbf{Conjecture}}

\newtheorem{hyp}{\textbf{Hypothesis}}

\newtheorem{lem}{\textbf{Lemma}}

\newtheorem{df}{\textbf{Definition}}

\newtheorem{cor}{\textbf{Corollary}}

\newtheorem{prop}{\textbf{Proposition}}

\newcommand{\A}{\Bbb{A}}

\newcommand{\C}{\Bbb{C}}

\newcommand{\HH}{\mathcal{H}}

\newcommand{\Q}{\Bbb{Q}}

\newcommand{\LL}{\mathcal{L}}

\newcommand{\II}{\mathcal{I}}

\newcommand{\FF}{\mathcal{F}}

\newcommand{\WW}{\mathcal{W}}

\newcommand{\TT}{\mathcal{T}}

\newcommand{\F}{\Bbb{F}}

\newcommand{\Z}{\Bbb{Z}}

\newcommand{\T}{\Bbb{T}}

\newcommand{\p}{\frak{p}}

\newcommand{\m}{\frak{m}}

\newcommand{\n}{\frak{n}}

\newcommand{\Wt}{\text{Wt}}

\newcommand{\Sp}{\text{Sp}}

\newcommand{\diag}{\text{diag}}

\newcommand{\Frob}{\text{Frob}}

\newcommand{\Ind}{\text{Ind}}

\newcommand{\ind}{\text{ind}}

\newcommand{\Hom}{\text{Hom}}

\newcommand{\cok}{\text{cok}}

\newcommand{\End}{\text{End}}

\newcommand{\la}{\langle}

\newcommand{\ra}{\rangle}

\newcommand{\im}{\text{im}}

\newcommand{\Gal}{\text{Gal}}

\newcommand{\GL}{\text{GL}}

\newcommand{\OO}{\mathcal{O}}

\newcommand{\y}{\hspace{6pt}}

\title{{\bf{The local Langlands correspondence in families and Ihara's lemma for $U(n)$}}}

\author{Claus M. Sorensen}

\begin{document}

\date{}

\maketitle

\begin{abstract} 
The goal of this paper is to reformulate the conjectural "Ihara lemma" for $U(n)$ in terms of the local Langlands correspondence in families $\tilde{\pi}_{\Sigma}(\cdot)$, as currently being developed by Emerton and Helm. The reformulation roughly takes the following form. Suppose we are given an irreducible mod $\ell$ Galois representation $\bar{r}=\bar{r}_{\m}$, which is modular of full level (and small weight), and a finite set of places $\Sigma$ -- none of which divide $\ell$. Then $\tilde{\pi}_{\Sigma}(r_{\m})$ exists, and has a global realization as a
natural module of algebraic modular forms, where $r_{\m}$ is the universal $\Sigma$-deformation of $\bar{r}$. This is unconditional for $n=2$, 
where Ihara's lemma is an almost trivial consequence of the strong approximation theorem.
\footnote{{\it{2000 AMS Mathematics Classification}}: 11F33.}
\end{abstract}

\onehalfspacing

\section{Introduction}

The modularity of elliptic curves $E_{/\Q}$ (and hence Fermat's Last Theorem) was proved by means of the Galois representation $r_{E,\ell}$ of $\Gamma_{\Q}=\Gal(\bar{\Q}/\Q)$ on the $\ell$-adic Tate module $T_{\ell}(E)\simeq \Z_{\ell}^{\oplus 2}$, for various primes $\ell$. The key step was to prove a modularity lifting theorem, roughly saying that 
$r_{E,\ell}$ arises from a modular form if its reduction $\bar{r}_{E,\ell}$ does. In a fancier language, one sought to identify a universal deformation ring $R_{\Sigma}^{univ}$ for 
$\bar{r}_{E,\ell}$ with a corresponding (localized) Hecke algebra $\T_{\Sigma}$ -- which parametrizes those deformations of $\bar{r}_{E,\ell}$ coming from modular forms. The set of primes $\Sigma$ is where one allows ramification. Wiles and Taylor found an ingenious method to prove $R_{\varnothing}^{univ}=\T_{\varnothing}$, by a certain intricate patching argument. This is the so-called minimal case. The general case, where $\Sigma\neq \varnothing$, is then reduced to the minimal case by verifying that both sides of the isomorphism $R_{\Sigma}^{univ}=\T_{\Sigma}$ grow in precisely the same way when one augments $\Sigma$. This reduction step requires one to control "level-raising" congruences on the Hecke side. The main ingredient which makes this work is a classical lemma of Ihara, which states that $\ker(J_0(N)^{\oplus 2}\rightarrow J_0(Np))$ 
is Eisenstein whenever $p \nmid N$.

\medskip

Clozel, Harris, and Taylor, almost succeeded in mimicking this approach for $\GL_n$ -- or, more accurately, for definite unitary groups in $n$ variables, say $G=U(B,\star)_{/F^+}$, which becomes an inner form $B^{\times}$ of $\GL_n$ over an imaginary CM-field $F=E F^+$ (assumed to be unramified everywhere over $F^+$). Fix a prime $\ell>n$, which splits in $E$, is unramified in $F^+$, and such that $S_{\ell}\cap S(B)=\varnothing$. Now, start off with an absolutely irreducible mod $\ell$ Galois representation,
$$
\bar{r}:\Gamma_F=\Gal(\bar{F}/F)\longrightarrow \GL_n(k).
$$
Here $k/\F_{\ell}$ is a finite extension, which we think of as the residue field $k=\OO/\lambda$ of a large enough coefficient field $K/\Q_{\ell}$. To make the modularity lifting machinery work, we need to assume $\bar{r}$ has large image -- more precisely, that $\bar{r}(\Gamma_{F^+(\zeta_{\ell})})$ is "big", in a formal sense (see Definition 2.5.1, p. 55 in \cite{1}). The fundamental hypothesis is of course that $\bar{r}$ is modular, by which we mean it "comes from" a modular form on $G$. That is, $\bar{r}\simeq \bar{r}_{\m}$ for a maximal ideal  
$\m\subset \T_{a,\{\rho_v\},\varnothing}^T(U)$ of the Hecke algebra acting faithfully on the module of algebraic modular forms, $S_{a,\{\rho_v\},\varnothing}(U,\OO)$.
We refer to the main text for any unexplained notation (see 3.2 and 3.3 below); here we will only point out that $\T^T$ is the Hecke algebra "away" from the finite set of split places of $F^+$,
$$
T=S \cup S_{\ell}\cup S(B), 
$$
for a collection of banal places $S$ (recall that $v$ is banal if $\#\GL_n(k(v))$ is not divisible by $\ell$). Moreover, for each $v \in T$ we choose a place $\tilde{v}$ of $F$ above it, which in particular pins down an identification $G(F_v^+)\simeq \GL_n(F_{\tilde{v}})$, well-defined up to conjugacy -- for $v \notin S(B)$. We must impose two rather strong conditions on these data (which can be traced back to \cite{1} -- linked to the modularity lifting techniques available at the time):

\begin{itemize}
\item $U=U(\varnothing) \subset G(\A_{F^+}^{\infty})$ is \underline{maximal} compact ("level one"),
\item $a=(a_{\tau,i})$ is in the "Fontaine-Laffaille" range ($0\leq a_{\tau,i}<\ell-n$).
\end{itemize}

\noindent There are additional mild hypotheses on the types $\{\rho_v\}$, which we suppress for now (see 4.2). The goal of  \cite{1} was to identify the universal deformation ring 
$R_{\Sigma}^{univ}$ of $\bar{r}$ with a certain localized Hecke algebra $\T_{\Sigma}$, following the strategy for $\GL_2$ outlined in the first paragraph. The minimal case ($\Sigma=\varnothing$) did not cause much trouble -- for the most part it is a notational challenge (the original approach goes through). Not surprisingly, when $\Sigma\subset S$ is enlarged, the main obstacle was to somehow measure the presence of congruences between modular forms of different levels $U(\Sigma)$ -- this group coincides with $U$ at all $v \notin \Sigma$, but at places $v \in \Sigma$ one replaces the maximal compact $U_v$ with the mirahoric subgroup $U_1(\tilde{v}^n)$ of depth $n$ -- the {\it{same}} $n$ as in $\GL_n$. We pull back $\m$ to a maximal ideal of level $U(\Sigma)$, and look at the universal modular $\Sigma$-deformation of $\bar{r}$,
$$
\text{$r_{\m}: \Gamma_F \longrightarrow \GL_n(\T_{\Sigma})$, $\y$ $\T_{\Sigma}=\T_{a,\{\rho_v\},\varnothing}^T(U(\Sigma))_{\m}$.}
$$
Here $\T_{\Sigma}$ is a {\it{reduced}} complete local Noetherian $\OO$-algebra. The main result of \cite{1}, which involves ideas from Mann's Harvard thesis \cite{11}, is an
isomorphism $R_{\Sigma}^{univ}=\T_{\Sigma}$ contingent on a conjectural (for $n>2$) analogue of Ihara's lemma.  Namely, Conjecture B in \cite{1}:

\begin{cn}
Let $\m \subset \T_{a,\{\rho_v\},\varnothing}^T(U)$ be a non-Eisenstein maximal ideal (that is, for which $\bar{r}_{\m}$ is absolutely irreducible). Let $x \in S$ be a banal place of $F^+$. Then every irreducible $\bar{k}[G(F_x^+)]$-submodule of $S_{a,\{\rho_v\},\varnothing}(U^x,\bar{k})_{\m}$ is generic.
\end{cn}

For $n=2$ this is almost trivially true, since non-generic representations of $\GL_2$ are characters, but for $n>2$ it is an important open problem -- whose arithmetic potential is not limited to $R=\T$ theorems, but conceivably has applications to special values of adjoint $L$-functions. It should be emphasized that Taylor has since shown that at least the reduced quotient $(R_{\Sigma}^{univ})^{red}$ is isomorphic to $\T_{\Sigma}$, by adapting Kisin's variant of the patching argument, which somehow does not distinguish between the minimal case and non-minimal case -- but treats both at once. However, this is a priori a weaker statement.

\medskip

The point of this paper is to revisit Ihara's lemma for $G$, and recast it in a more robust contemporary framework -- invoking the local Langlands correspondence in families, and its mod $\ell$ variant, as developed recently by Emerton and Helm. We briefly give the gist of it. The question is whether one can interpolate the local Langlands correspondence in families, and somehow make sense of it over more general coefficient rings $A$ -- as opposed to just fields of characteristic zero. To fix ideas, we take $A$ to be a complete {\it{reduced}} Noetherian local $\OO$-algebra. Suppose we have a collection of Galois representations $r_w: \Gamma_{F_w} \rightarrow \GL_n(A)$; one for each place $w \in \Delta$, where $\Delta$ is a finite set of finite places of a number field $F$ -- none of which lie above the residue characteristic $\ell$. For each minimal prime $\p \subset A$, one can look at the specializations $r_w\otimes_A \kappa(\p)$, over the residue field $\kappa(\p)$ (which is a finite extension of $K$), and via local Langlands attach a smooth $\kappa(\p)$-representation $\tilde{\pi}(r_w\otimes_A \kappa(\p))$ of $\GL_n(F_w)$. Examples show that it is much more natural to take $\pi(\cdot)$ to be the "generic" local Langlands correspondence of Breuil and Schneider -- thus the representations are always generic, but possibly reducible; the tilde in $\tilde{\pi}(\cdot)$ signifies taking the smooth $\kappa(\p)$-dual -- a convenient normalization. Now, it is natural to ask for a smooth $A[\GL_n(F_{\Delta})]$-module $H$ which {\it{interpolates}} these representations, as $\p$ varies (in the sense of (2) below). 
Theorem 6.2.1 in \cite{5} shows that there can be at {\it{most}} one such $H$ satisfying two additional technical hypotheses -- (1) and (3) below,

\begin{itemize}
\item[(1)] $H$ is $A$-torsion free (only zerodivisors can annihilate a nonzero $h \in H$).
\item[(2)] For each minimal prime $\p \subset A$, there is a $\GL_n(F_{\Delta})$-equivariant isomorphism,
$$
\otimes_{w\in \Delta} \tilde{\pi}(r_w\otimes_A \kappa(\p)) \overset{\sim}{\longrightarrow} \kappa(\p)\otimes_A H.
$$
\item[(3)] The smooth dual $(H/\m H)^{\vee}$ is {\it{essentially AIG}} (absolutely irreducible and generic) -- which means its socle is absolutely irreducible, generic, and all other constituents are non-generic.
\end{itemize}

\noindent If such an $H$ exists, one denotes it $\tilde{\pi}(\{r_w\}_{w\in \Delta})$. In on-going work of Helm, he shows the existence of $H$ when $\Delta$ consists of banal places -- by linking deformation rings to his theory of the integral Bernstein center. See Theorem 7.8 and Example 7.10 of \cite{15}. 

\medskip

In our setup we will take $A=\T_{\Sigma}$, and $\Delta=\tilde{\Sigma}$ (the collection of all choices $\tilde{v}$, for $v \in \Sigma$). We then seek a global realization of 
the module $\tilde{\pi}(\{r_{\m,\tilde{v}}\}_{v\in \Sigma})$. Our main result exhibits it as a module of algebraic modular forms, granting Conjecture 1 (and a minor multiplicity one assumption):

\begin{thm}
Admit Ihara's lemma, and multiplicity one for the automorphic spectrum of $G$ -- both hold for $n=2$.
Then there are $G(F_{\Sigma}^+)$-equivariant $\T_{\Sigma}$-linear isomorphisms, unique modulo $\T_{\Sigma}^{\times}$ by 6.2.1 (5) in \cite{5},
$$
\tilde{\pi}(\{r_{\m,\tilde{v}}\}_{v \in \Sigma})\simeq (\otimes_{v \in \Sigma}\tilde{\pi}(r_{\m,\tilde{v}}))^{tf}\simeq S_{a,\{\rho_v\},\varnothing}(U^{\Sigma},\OO)_{\m},
$$
where $tf$ denotes the maximal $\T_{\Sigma}$-torsionfree quotient (and $U^{\Sigma}=\prod_{v\notin \Sigma}U_v$ is a product of hyperspecial maximal compact subgroups, away from $\Sigma$). 
\end{thm}

The first isomorphism in the Theorem is Proposition 6.2.4 in \cite{5}. The point of our paper is to give a concrete global realization of $\tilde{\pi}(\{r_{\m,\tilde{v}}\}_{v \in \Sigma})$, contingent on Ihara's lemma (and multiplicity one -- which should be a minor issue in comparison). In fact, the second isomorphism in Theorem 1 is {\it{equivalent}} to Ihara's lemma, as we will now explain.

\medskip

From now on, for brevity, we will employ the notation
$$
H=S_{a,\{\rho_v\},\varnothing}(U^{\Sigma},\OO)_{\m}.
$$
It becomes a $\T_{\Sigma}$-module in a natural way (see 5.2). Thus the content of Theorem 1 is that $H$ satisfies the desiderata (1)--(3) above. We first show that $H^{tf}$ satisfies these criteria, and only at the very end we show that in fact $H=H^{tf}$ is torsion-free (over $\T_{\Sigma}$). It is fairly easy to show $H$ satisfies (2), and deduce that $H^{tf}$ does too. What requires some more serious thought is to show that $(H/\m H)^{\vee}$ is essentially AIG -- and this is where Ihara's lemma comes in. Indeed, it {\it{implies}} Ihara's lemma: Petersson duality identifies $(H/\m H)^{\vee}$ with a space of mod $\ell$ modular forms on a closely related unitary group $G'=U(B^{op},\star)$ (see Theorem 4) --
and all simple submodules hereof must be generic, because the socle is. Conversely, we must show $(H/\m H)^{\vee}$ has a {\it{unique}} generic constituent -- which must sit as a submodule if we believe Ihara. To do that, we introduce a certain (exact) localization functor $\LL=\LL_{n,\Sigma}$, which commutes with base change: Take $U_1(\tilde{v}^n)$-invariants, for all $v \in \Sigma$, and localize at the ideal generated by Atkin-Lehner operators $U^{(j)}$. We study its local properties in the first section, building on Mann's thesis
\cite{11}; with focus on how it detects generic representations. Applying $\LL$ to $H/\m H$, we essentially reduce to verifying $\LL(H)$ is free of rank one over $\T_{\Sigma}$.
The freeness is a by-product of the proof of $R_{\Sigma}^{univ}\overset{\sim}{\longrightarrow} \T_{\Sigma}$, along the lines of \cite{2} -- by a {\it{second}} application of Ihara's lemma. This freeness is what used to be called mod $\ell$ "multiplicity one". The fact that the rank is one obviously uses the $m_{\pi}=1$ hypothesis for $G$, in conjunction with our local results on $\LL$ combined with known properties of the mod $\ell$ local Langlands correspondence of Emerton and Helm, \cite{5}.

\medskip

In the last Chapter of this paper we give a mod $\ell$ reformulation of Ihara's lemma, in terms of the mod $\ell$ local Langlands correspondence $\bar{\pi}(\cdot)$ of Emerton-Helm. As a Corollary of Theorem 1, we show that (admitting Ihara) $\otimes_{v\in \Sigma} \tilde{\bar{\pi}}(\bar{r}_{\m,\tilde{v}})$ maps onto $H/\m H$. See Corollary 4 in section 6 for more details.

\medskip

We anticipate that Theorem 1 will have applications to strong local-global compatibility for unitary groups $G$, in two variables, in the vein of \cite{9}, which is concerned with the case of $\GL(2)_{/\Q}$. The key object is the completed cohomology $\hat{H}^0(U^{\Sigma\cup S_{\ell}})_{\OO,\m}$, where as indicated we shrink the level $U_{\ell}\rightarrow 1$.
It becomes a module for the "big" Hecke algebra $\T_{\Sigma}^{big}$; the corresponding inverse limit over shrinking $U_{\ell}$. The goal is to somehow factor the two natural actions of 
$G(F_{\Sigma}^+)$ and $G(F^+\otimes_{\Q}\Q_{\ell})$ as a (completed) tensor product,
$$
\hat{H}^0(U^{\Sigma\cup S_{\ell}})_{\OO,\m} \simeq \pi(\{r_{\m,\tilde{v}}\}_{v\in S_{\ell}})
\overset{\curlywedge}{\otimes}_{\T_{\Sigma}^{big}} \pi(\{r_{\m,\tilde{v}}\}_{v\in \Sigma}).
$$
(Here $\overset{\curlywedge}{\otimes}$ is a direct limit of $\varpi$-adically completed tensor products; see Definition C. 43 on p. 114 in \cite{9}.)
The first factor $\pi(\{r_{\m,\tilde{v}}\}_{v\in S_{\ell}})$ arises from the "$\ell$"-adic local Langlands correspondence for $\GL_2(\Q_{\ell})$ -- of Berger, Breuil, Colmez, Kisin, Paskunas and others. The other factor $\pi(\{r_{\m,\tilde{v}}\}_{v\in \Sigma})$ should be the smooth $\OO$-{\it{dual}} of local Langlands in families over  $\T_{\Sigma}^{big}$ --
a "codamissible" module in the terminology of Emerton's Appendix C of \cite{9}. To show this, one would introduce the module
$$
X=\Hom_{\T_{\Sigma}^{big}[G(F^+\otimes_{\Q}\Q_{\ell})]}(\pi(\{r_{\m,\tilde{v}}\}_{v\in S_{\ell}}), \hat{H}^0(U^{\Sigma\cup S_{\ell}})_{\OO,\m})
$$
and try to identify its maximal cotorsion-free submodule $X_{ctf}$ with $\pi(\{r_{\m,\tilde{v}}\}_{v\in \Sigma})$. Among other things, the key property one would have to verify is that
the representation $(X/\varpi X)[\m]$ is essentially AIG. However, the latter embeds into
$$
\Hom_{k[G(\OO_{F^+}\otimes_{\Z}\Z_{\ell})]}(W_a,H^0(U^{\Sigma\cup S_{\ell}})_{k,\m}[\m])\simeq S_{a,\varnothing}(U^{\Sigma},k)_{k,\m}[\m]
$$
for suitable Serre weights $a$ -- and our paper shows this last representation on mod $\ell$ modular forms is indeed essentially AIG (hence so is any submodule). We hope to work out the details in a separate "companion" paper with P. Chojecki, in natural continuation of \cite{10}. A preprint is now available \cite{16}, in which we prove strong local-global compatibility when $\bar{r}$ is irreducible at all places above $\ell$ -- admittedly a rather strong assumption, which we expect can be weakened significantly.

\medskip

On that note, some of the initial motivation for writing this paper was to ease the transition between \cite{9} and \cite{5}: In part (1) of Theorem 4.1.2 on p. 41 of \cite{9}, 
the Kirillov functor $\FF_{\infty}$ (in the notation of 2.4 below) is claimed to be exact, without proof. It turns out, for $n=2$ it {\it{is}} indeed exact -- see Corollary 2 in section 2.5 below -- but this is not immediate, and it is not obvious at all if exactness of $\FF_{\infty}$ continues to hold for $n>2$ (left-exactness is easy, though). What's worse, it is unclear to what extent it commutes with extension of scalars. This leads to an ad hoc notion of "generic" in Definition 4.1.3 of \cite{9} -- used in his Theorems 4.4.1 and 5.7.7, which a priori could be different from the standard notion of generic (existence of Whittaker models) -- used in \cite{5}, and therefore the connection to local Langlands in families is somewhat blurry (at least to the untrained eye) in the current version of \cite{9} (as of April, 2014). One of the points of this paper is to replace $\FF_{\infty}$ by a certain localization functor $\LL_n$, which enjoys all the desired properties. To verify these we rely on Mann's Harvard thesis \cite{11}, and extend some of Mann's results.

\medskip

\noindent {\it{Acknowledgements}}. I wish to thank P. Chojecki for useful correspondence, which eventually led to the writing of this paper, and for his continued interest.
I also wish to acknowledge the impact of many enjoyable conversations with D. Helm, F. Herzig, S. W. Shin, and L. Xiao on related topics. Finally, I thank M. Emerton for  clarifications. 

\section{Local preliminaries}

Throughout this section, we work with two local fields. The base field is a finite extension $F/\Q_p$, with integers $\OO_F$, uniformizer $\varpi=\varpi_F$, and residue field
$k_F$. The coefficient field is a finite extension $K/\Q_{\ell}$, with integers $\OO$, maximal ideal $\lambda=(\varpi_K)$, and residue field $k$. Here $p \neq \ell$. Later on we will always assume $K$ is large enough to contain the image of every embedding $F \hookrightarrow \bar{K}$. We consider smooth representations of $G=\GL_n(F)$ on 
$A$-modules, where $A$ is an $\OO$-algebra. Eventually we will take $A$ to be a complete local Noetherian ring (specifically a localized Hecke algebra), but we will not need this assumption in the beginning. We will always assume we are in the \underline{banal} case: $\#\GL_n(k_F)$ is prime to $\ell$.

\subsection{Smooth representations in families}

Let $V$ be a smooth $A[G]$-module; which means $V$ is the union of all $U$-invariants $V^U$, where $U$ ranges over the compact open subgroups of $G$. Each $V^U$
is an $A$-submodule, with an action of the Hecke algebra $\HH_U$ of compactly supported $U$-biinvariant functions $G \rightarrow \OO$ (equipped with convolution). As an
$\OO$-module, $\HH_U$ is generated by characteristic functions of double cosets, $[U\alpha U]$ with $\alpha \in G$. The inclusion $V^U \subset V$ has a natural splitting: Define
an operator $e_U$ on $V$ by 
$$
e_U(x)=\frac{1}{[U:U_x]}\cdot {\sum}_{u \in U/U_x}u x,
$$
where $U_x \subset U$ is an open subgroup fixing $x$. This is well-defined since the pro-order $|U|$ is prime to $\ell$, and therefore the index $[U:U_x]$ is invertible in $\OO$.
Clearly this defines a retraction $e_U: V \rightarrow V^U$, and therefore $V=V^U \oplus \ker(e_U)$ as $A$-modules.

\begin{lem}
The functor $V \mapsto V^U$ is exact, and commutes with base change.
\end{lem}

\noindent {\it{Proof}}. Exactness is clear (make use of $e_U$). Now, suppose $B$ is an $A$-algebra, and we extend scalars to $B$. On one hand,
$$
V=V^U \oplus \ker(e_U) \Longrightarrow V\otimes_A B= (V^U \otimes_A B) \oplus (\ker(e_U)\otimes_A B).
$$
On the other hand, letting $e_U$ act on $V\otimes_A B$, we have 
$$
V\otimes_A B=(V\otimes_A B)^U \oplus \ker(V\otimes_A B \overset{e_U}{\longrightarrow} V\otimes_A B).
$$
There are obvious inclusions among the summands in these decompositions. Those inclusions must be equalities. $\square$

\subsection{Mirahoric subgroups and their invariants}

We will now specialize $U$ to be one of the mirahoric subgroups $U_1(\varpi^r)$, which appear in the definition of conductors, for example. Recall their definition,
$$
U_1(\varpi^r)=\{\text{$u \in \GL_n(\OO_F)$: $u\equiv \begin{pmatrix} * & * \\ 0 & 1\end{pmatrix}$ (mod $\varpi^r$) }\}.
$$
More precisely, the last row of $u$ is congruent to $(0,\ldots,0,1)$ modulo $\varpi^r$. We extend the definition to $r=\infty$ by taking $U_1(\varpi^{\infty})$ to be the mirabolic subgroup $P(\OO_F)$, the intersection of all the $U_1(\varpi^r)$. For later use, observe that $U_1(\varpi^r)$ is generated by $P(\OO_F)$ and the principal congruence subgroup
$U(\varpi^r)$ consisting of $u\equiv 1$ (mod $\varpi^r$). For $V$ as above, we get an increasing sequence of $A$-modules,
$$
V^{U_1(\varpi)} \subset V^{U_1(\varpi^2)} \subset \cdots \subset V^{U_1(\varpi^r)}\subset \cdots \subset V^{P(\OO_F)}={\bigcup}_{r=1}^{\infty} V^{U_1(\varpi^r)}.
$$
Here $V^{U_1(\varpi^r)}$ has an $A$-linear action of the Hecke algebra $\HH_{U_1(\varpi^r)}$. We specify certain Hecke operators $U^{(j)}$ in here, for $j=1,\ldots,n-1$, analogous to the Atkin-Lehner operator $U_p$ on classical elliptic modular forms of level $\Gamma_0(p)$.

\begin{df}
$U^{(j)}=[U_1(\varpi^r) \alpha_j U_1(\varpi^r)]$, where $\alpha_j=\begin{pmatrix} \varpi \cdot 1_j & \\ & 1_{n-j}\end{pmatrix}$.
\end{df}

\noindent Their action can be made explicit, as was done in Mann's Harvard PhD thesis (among other things). Unfortunately, this is neither published nor widely available. 
(The author is in possession of a copy.) Many parts of it are reviewed or reproved in \cite{1}. One of Mann's key findings is a set of coset representatives, which is independent of $r$. We recall this set here for future reference.

\begin{prop}
Let $r>0$ and $j<n$. Then there is a partition
$$
U_1(\varpi^r) \alpha_j U_1(\varpi^r)=\bigsqcup_{I,b} bU_1(\varpi^r),
$$
where $I \subset \{1,\ldots,n-1\}$ is a subset of cardinality $|I|=j$, and $b \in B$ ranges over the set $B_I$ of upper-triangular matrices in $G$ satisfying the conditions below:
\begin{itemize}
\item $b_{ii}=\varpi$ for $i \in I$, and $b_{ii}=1$ otherwise.
\item For a pair $i<j$ such that $i\in I$ and $j \notin I$, allow arbitrary $b_{ij} \in Z$.
Also, $b_{ij}=0$ for any other pair $i<j$. (That is, such that $i \notin I$ or $j \in I$.)
\end{itemize}
Here $Z \subset \OO_F$ is a complete set of representatives for $k_F$, containing $0$.
\end{prop}

\noindent {\it{Proof}}. This is Proposition 4.1 on p. 10 in \cite{11}. See also p. 82 in \cite{1}. $\square$

\medskip

\noindent We allow $r=\infty$ in the Proposition. For $r=0$ Mann gives a similar set of representatives; the only difference is that $I$ runs over subsets of
$\{1,\ldots,n\}$. Using this, he deduces the following from the fact that $\HH_{\GL_n(\OO_F)}$ is commutative: 

\begin{cor}
The operators $U^{(j)}$, for $j=1,\ldots,n-1$, generate a commutative subalgebra of $\HH_{U_1(\varpi^r)}$.
\end{cor}

\noindent {\it{Proof}}. This is Corollary 4.2 on p. 11 in [Man]. $\square$

\medskip

\noindent This subalgebra is a quotient of the polynomial algebra $\OO[X_1,\ldots,X_{n-1}]$, which we will just denote $\OO[X]$ in what follows. In the notation above, each $V^{U_1(\varpi^r)}$ thus becomes a module over $A[X]=\OO[X]\otimes_{\OO}A$, by sending $X_j \mapsto U^{(j)}$.

\begin{lem}
For $0<r\leq s$, the inclusion $V^{U_1(\varpi^r)}\subset V^{U_1(\varpi^s)}$ is $A[X]$-linear. The corresponding quotient module is $A$-torsion free if $V$ is.
\end{lem}

\noindent {\it{Proof}}. The inclusion is clearly $A[X]$-linear by Mann's Proposition (the independence of $r$). The torsion freeness is easily checked: If $ax$ is fixed by $U_1(\varpi^r)$,
we must have $a(ux-x)=0$ for all $u \in U_1(\varpi^r)$. If $V$ is $A$-torsion free we deduce $ux=x$, unless $a$ is a zero-divisor. $\square$

\subsection{Localization at the Atkin-Lehner ideal}

We will view the various $A[X]$-modules $V^{U_1(\varpi^r)}$, from the previous section, simply as $\OO[X]$-modules (via $\OO\rightarrow A$) and localize them at the maximal ideal
$$
\n=(\lambda,X_1,\ldots,X_{n-1})\subset \OO[X].
$$
(Recall that $\lambda \subset \OO$ is the maximal ideal, and $k=\OO/\lambda$.) We shall refer to $\n$ as the Atkin-Lehner ideal, again by analogy with the classical theory.
More formally, we introduce a sequence of functors $\LL_r$ from smooth $A[G]$-modules to the category of $A$-modules with $\OO[X]_{\n}$-action, as follows:

\begin{df}
$\LL_r(V)=(V^{U_1(\varpi^r)})_{\n}=V^{U_1(\varpi^r)} \otimes_{\OO[X]}\OO[X]_{\n}$.
\end{df}

\noindent This is motivated by the definition of $X_S$ at the bottom of p. 147 in \cite{1}.
Let us summarize some of the main properties of $\LL_r$, which will be used below.

\begin{lem}
\begin{itemize}
\item[(1)] $\LL_r$ is exact, and commutes with base change.
\item[(2)] For $0<r \leq s$, the inclusion $\LL_r(V)\subset \LL_s(V)$ is $\OO[X]_{\n}$-linear, and the corresponding quotient module is $A$-torsion free if $V$ is.
\end{itemize}
\end{lem}

\noindent {\it{Proof}}. For (1), note that $\LL_r$ is a composition of two functors with those properties; taking $U_1(\varpi^r)$-invariants and then localizing at $\n$. Therefore,
$$
\LL_r(V\otimes_A B)=\OO[X]_{\n}\otimes_{\OO[X]} (V^{U_1(\varpi^r)}\otimes_A B)=\LL_r(V)\otimes_A B,
$$
and exactness is immediate. For (2), note that the quotient in question is 
$$
\LL_s(V)/\LL_r(V)\simeq (V^{U_1(\varpi^s)}/V^{U_1(\varpi^r)})_{\n}=M_{\n},
$$
and $M$ is known to be $A$-torsion free, if $V$ is. Then so is $M_{\n}$, as is easily verified. (Indeed, if $\frac{m}{s}\in M_{\n}$
is annihilated by $a \in A$, then $s'am=0$, for some $s' \notin \n$. Hence $a|0$ or $s'm=0$. The latter implies $\frac{m}{s}=0$ in $M_{\n}$.)
$\square$

\bigskip

Over $\bar{K}$ there is the following more concrete description: If $V$ is an admissible $\OO[G]$-module (which means the $U$-invariants are finitely generated over $\OO$), and $K$ is large enough, $\LL_r(V)\otimes_{\OO}\bar{K}$ breaks up into
generalized eigenspaces in $V^{U_1(\varpi^r)}\otimes_{\OO}\bar{K}$ for eigensystems $\OO[X]\rightarrow \OO$ sending each $X_j$ to an eigenvalue in $\lambda\subset \OO$. 
(That is, systems lifting the natural mod $\n$ reduction map $\OO[X]\rightarrow k$.)

\subsection{The common kernel of the $U$-operators}

We introduce variants of the functors $\LL_r$, which are easier to compute, but they may not share the same nice properties (they turn out to be left-exact, but in general, for $n>2$, it is not clear whether they are exact, or whether they commute with base change). 

\begin{df}
$\FF_r(V)=\{x \in V^{U_1(\varpi^r)}: U^{(1)}x=\cdots=U^{(n-1)}x=0\}$.
\end{df}

\noindent This defines a functor from smooth $A[G]$-modules to $A$-modules. Again, we include the case $r=\infty$, where one takes the annihilator of all the $U^{(j)}$ in
$V^{P(\OO_F)}$. We note that $\FF_{\infty}$ generalizes what Emerton calls the Kirillov functor (for $n=2$) in section 4.1 of [Eme]. Clearly, as $r$ varies, we have a chain 
$$
\FF_1(V)\subset \FF_2(V)\subset \cdots \subset \FF_r(V) \subset \cdots \subset \FF_{\infty}(V)={\bigcup}_{r=1}^{\infty}\FF_r(V),
$$
and the subquotients are $A$-torsion free if $V$ is. 

\begin{lem}
The canonical map $\FF_r(V) \rightarrow \LL_r(V)$ is injective.
\end{lem}

\noindent {\it{Proof}}. Suppose $x \in \FF_r(V)$ maps to $0$ in $\LL_r(V)$. That is, $sx=0$ for some $s \notin \n$. Since $x$ is annihilated by every $U^{(j)}$, the algebra $\OO[X]$ acts on $x$ via the constant-term map $\gamma:\OO[X]\rightarrow \OO$ (which sends $X_j \mapsto 0$). In particular, $\gamma(s)x=0$; from which we deduce that $x=0$ or $\gamma(s) \in \lambda$ (that is, is a non-unit in $\OO$). The latter violates the fact that $s \notin \n$. $\square$

\medskip

At the end of the next subsection (Corollary 3) we give an example showing that the inclusion $\FF_r(V) \hookrightarrow \LL_r(V)$ is usually {\it{not}} surjective: If $V$ is an unramified principal series for $\GL(2)$ over $\bar{K}$, the eigenspace $\FF_r(V)$ is one-dimensional, whereas the {\it{generalized}} eigenspace $\LL_r(V)$ has dimension $r-1$.

\subsection{A refinement of a result of Mann}

To start off, let us get our definition of {\it{generic}} in place. Here we follow \cite{17}.

\begin{df}
A (possibly reducible) smooth admissible $G$-representation $V$ is said to be generic if there is a $G$-equivariant embedding
$V \hookrightarrow \Ind_N^G(\psi)$, for some (hence any) non-trivial additive character $\psi$ of $F$.
\end{df}

The main result of \cite{17} (that is, Proposition 3.2 on p. 114) shows that parabolically induced representations $V=\Ind_Q^G(\tau_i)$ are generic if the $\tau_i$ are ordered such that
they have decreasing exponents (representations of "Langlands type"). In fact, \cite{17} proves the existence of a Kirillov model -- that is, for Whittaker functions $W$ in the image of $V \hookrightarrow \Ind_N^G(\psi)$, the mirabolic restriction map $W \mapsto W|_P$ is {\it{injective}}. (Their proof is over $\C$, but works well in banal characteristic.)
This was announced as Theorem 5.20, p. 50 in \cite{6}, for irreducible $V$ (and proved in a sequel) -- which was a conjecture of Gelfand and Kazhdan. 

\medskip

The following will be the crucial local input used later in our paper.

\begin{thm}
Allow $A=\bar{K}$ or $A=\bar{k}$ in part (1), and let $A=\bar{k}$ in parts (2)--(4):
\begin{itemize}
\item[(1)] Let $V$ be a generic representation of $G$ over $A$, possibly reducible. Then $\FF_{r}(V)=\FF_{\infty}(V)$ is one-dimensional over $A$, for $r$ sufficiently large.
\item[(2)] (Mann) Suppose $V=i_B(\bar{\chi}_i)$ is an unramified principal series representation of $G$ over $\bar{k}$, possibly reducible. Then $\FF_n(V)=\LL_n(V)$ is one-dimensional over $\bar{k}$ (where the subscript $n$ is the \underline{same} $n$ as in $\GL_n$). 
\item[(3)] Suppose $V=i_B(\bar{\chi}_i)$ is an unramified principal series representation of $G$ over $\bar{k}$, and let $V_+$ be its unique generic constituent. Then 
$\FF_n(V_+)=\LL_n(V_+)$ is one-dimensional over $\bar{k}$.
\item[(4)] Suppose $V=i_B(\bar{\chi}_i)$ is an unramified principal series representation of $G$ over $\bar{k}$, and let $V'$ be a non-generic constituent. Then 
$\FF_n(V')=\LL_n(V')$ is trivial.

\end{itemize}
\end{thm}

\noindent {\it{Proof}}. Part (1) is a computation in the Kirillov model. Pick a non-trivial additive character $\psi: F \rightarrow A^{\times}$, with kernel $\OO_F$, and view it as a character of the unipotent radical $N \subset B$ in the usual fashion (via $N^{ab}\simeq F^{n-1}$ composed with summation). Let $\Ind_N^G(\psi)$ be its smooth induction to $G$,
and $\ind_N^G(\psi)$ its compact induction. By our genericity assumption, $V \hookrightarrow \Ind_N^G(\psi)$, but $V$ need not be irreducible. Thus we identify $V$ with a space of Whittaker functions, $\WW_V$. We restrict each $W \in \WW_V$ to the mirabolic subgroup $P=P(F)$. As discussed in the paragraph leading up to the Theorem, the map $W \mapsto W|_P$ is injective -- by the result of Bernstein-Zelevinsky: Indeed, if there was a nonzero kernel in $\WW_V$, it would contain a simple submodule since $V$ has finite length. Apply Bernstein-Zelevinsky to this sub, which is necessarily generic.   

\medskip

So,
$V \hookrightarrow \Ind_N^P(\psi)$. Essentially by duality, $V|_P$ also contains a copy of the (irreducible) Kirillov representation, $\ind_N^P(\psi)\hookrightarrow V$ (see p. 49 in \cite{6}; look at the irreducible generic sub $V_+\subset V$). We will prove (1) by verifying that $\FF_{\infty}(\Ind_N^P(\psi))$ is spanned by one nonzero function $W_{\infty}$, which is compactly supported modulo $N$. This proves (1) since any $P$-map $\ind_N^P(\psi) \hookrightarrow \Ind_N^P(\psi)$ is proportional to the standard embedding (this is part (3) of the Proposition on p. 49 in \cite{6}).

\bigskip

With every integer tuple $m=(m_1,\ldots,m_{n-1})\in \Z^{n-1}$, we associate the diagonal matrix
$$
\rho_m=\diag(\varpi^{m_1},\ldots, \varpi^{m_{n-1}},1)\in P.
$$
Using the Iwasawa decomposition for $\GL_{n-1}$, it is not hard to check that these $\rho_m$ form a complete set of inequivalent representatives for the double cosets
$$
N\backslash P/P(\OO_F) \leftrightarrow \{\rho_m\}.
$$
Consequently, every $P(\OO_F)$-invariant $W \in \Ind_N^P(\psi)$ is completely determined by its values $W(\rho_m)$. Furthermore, as an easy calculation shows,
$$
W(\rho_m)\neq 0 \Longrightarrow \psi(N \cap \rho_m P(\OO_F) \rho_m^{-1})=1 \Longrightarrow m_1 \geq \cdots \geq m_{n-1}\geq 0.
$$
That is, $m \in \Z_+^{n-1}$ must be dominant (here we use that $\psi$ was chosen to have conductor $\OO_F$). This sets up an isomorphism of $A$-vector spaces,
$$
\text{$\Ind_N^P(\psi)^{P(\OO_F)} \overset{\sim}{\longrightarrow} A^{\Z_+^{n-1}}$, $\y$ $W \mapsto (W(\rho_m))_{m \in \Z_+^{n-1}}$,}
$$
which identifies $\ind_N^P(\psi)^{P(\OO_F)}$ with finitely supported sequences, $\oplus_{\Z_+^{n-1}}A$. 

We will transfer the $U^{(j)}$-action across this isomorphism, and explicitly write down the corresponding operator on tuples. Say $a=(a_m)\in A^{\Z_+^{n-1}}$. Then,
using Mann's double coset representatives given in Proposition 1, 
$$
(U^{(j)}a)_m={\sum}_{I\subset \{1,\ldots,n-1\}: |I|=j} |B_I|\cdot a_{m+e_I},
$$
where $a_{m+e_I}=0$ if $m+e_I$ is not dominant. In this expansion $|B_I|$ is a non-negative power of $q_F$, and $e_I$ is the characteristic function of $I$.
It follows from Proposition 1, with $r=\infty$: Indeed, for any $P(\OO_F)$-invariant $W \in \Ind_N^P(\psi)$, 
$$
(U^{(j)}W)(\rho_m)=\sum_{I,b}W(\rho_m b)=\sum_{I,b}W(\rho_m b \rho_m^{-1}\rho_m),
$$
and we factor $\rho_m b\rho_m^{-1}=ut$, where $u \in N$, and $t=\diag(b_{11},\ldots,b_{nn})\in T$. Then
$t \rho_m=\rho_{m+e_I}$, where $e_I$ is the characteristic function of $I$. Moreover, we have
$$
\psi(u)=\psi(\varpi^{m_1-m_2}b_{12}b_{22}^{-1}+\cdots+\varpi^{m_{n-1}-0}b_{n-1,n}b_{nn}^{-1}).
$$
By the second bullet in Proposition 1, $b_{i-1,i}=0$ unless $i \notin I$ and $i-1 \in I$. In particular,
$b_{i-1,i}\neq 0$ implies $b_{ii}=1$. Since $m$ is dominant, the sum lies in $\OO_F$, so that $\psi(u)=1$. From this, it is straightforward to deduce the expression for 
$U^{(j)}a$. We proceed to compute the common kernel of the $U^{(j)}$-operators on tuples. We start with the last operator, $U^{(n-1)}$. Here $I$ would have to be the whole set
$\{1,\ldots,n-1\}$, and 
$$
(U^{(n-1)}a)_m=|B_{\{1,\ldots,n-1\}}| \cdot a_{m+(1,\ldots,1)}.
$$
We read off that $\ker U^{(n-1)}$ consists of $a$ such that $a_m=0$ whenever $m_{n-1}>0$. Move on to $U^{(n-2)}$, and intersect its kernel with $\ker U^{(n-1)}$.
For a tuple $a$ in the latter, only $I=\{1,\ldots,n-2\}$ contributes to the sum defining $U^{(n-2)}a$, so
$$
(U^{(n-2)}a)_m=|B_{\{1,\ldots,n-2\}}| \cdot a_{m+(1,\ldots,1,0)}.
$$
Again, we read off that $\ker U^{(n-2)}\cap \ker U^{(n-1)}$ consists of $a$ such that $a_m=0$ whenever $m_{n-1}>0$ or $m_{n-2}>m_{n-1}$. Continuing this way, eventually we find that the common kernel of all the $U^{(j)}$, as $j=1,\ldots,n-1$, consists of $a$ such that $a_m=0$ whenever at least one of the inequalities $m_{i-1}\geq m_i$ is strict (with the convention $m_n=0$). In other words, $m \neq 0$. Altogether, this shows that $\FF_{\infty}(\Ind_N^P(\psi))$ is spanned by the function $W_{\infty}$ with $W_{\infty}(\rho_m)=\delta_{m,0}$, which clearly is compactly supported modulo $N$, and therefore lies in $\ind_N^P(\psi)$, as desired. We take $r$ large enough such that $W_{\infty}$ lies in
$\FF_r\subset \FF_{\infty}$.

\bigskip

Part (2) is due to Mann: Since $A=\bar{k}$, we can think of $\LL_n(V)$ as just being the simultaneous generalized eigenspace for the $U^{(j)}$-operators on $V^{U_1(\varpi^n)}$, with eigenvalues $0$. By Proposition 4.4 on p. 12 in \cite{11}, these $e$-spaces are parametrized by proper subsets $\TT \subset \{1,\ldots,n\}$. The corresponding eigenvalue of $U^{(|\TT|)}$ is nonzero unless $\TT=\varnothing$ (by the explicit eigenvalue formula in loc. cit.). Thus $\LL_n(V)$ corresponds to $\TT=\varnothing$, and therefore has dimension 
one -- again, using the dimension formula in Mann's Proposition 4.4. (See also the proof of Corollary 4.5 in \cite{11}.) In particular, being one-dimensional, the generalized eigenspace $\LL_n(V)$ coincides with the actual eigenspace $\FF_n(V)$.

\bigskip

For part (3), first observe that since $\LL_n$ is exact, $\dim \LL_n(V_+)$ can be at most one -- using that $\LL_n(V)$ is one-dimensional, as explained in the previous paragraph. Moreover, by Lemma 4,
$$
\dim \FF_n(V_+) \leq \dim \LL_n(V_+)\leq 1.
$$
It remains to show $\FF_n(V_+)$ is nonzero: By permuting the $\bar{\chi}_i$, we may assume that $V_+\subset V=i_B(\bar{\chi}_i)$, without loss of generality. By part (1), we know
$\FF_{\infty}(V_+)$ is one-dimensional. Say, spanned by $\xi$. In addition, combining (1) and (2), we know that $\FF_{\infty}(V)=\FF_n(V)$. Indeed, $V=i_B(\bar{\chi}_i)$ is a (possibly reducible) generic representation. Thus, $\xi$ lies
in $\FF_n(V)$, and is therefore $U_1(\varpi^n)$-invariant. A fortiori, we have a nonzero vector $\xi \in \FF_n(V_+)$ as claimed. 
 
\bigskip

Part (4) now follows more or less immediately from the exactness of $\LL_n$. Indeed, again from (2) we get the estimate
$$
1=\dim \LL_n(V)\geq \dim\LL_n(V_+)+\dim \LL_n(V'),
$$
and we know $\dim\LL_n(V_+)=1$ from part (3). We conclude that $\dim \LL_n(V')=0$. $\square$

\begin{cor}
Let $A=\bar{K}$ or $A=\bar{k}$, as in the previous Theorem. Let $V$ be any (smooth) representation of $G$ over $A$. Then,
$$
\dim \FF_{\infty}(V)\geq \dim V_{N,\psi},
$$ 
where $V_{N,\psi}$ denotes the $(N,\psi)$-coinvariants of $V$, with $\psi$ an arbitrary generic character of $N^{ab}\simeq F^{n-1}$ (the Bernstein-Zelevinsky "top derivative" $V^{(n)}$).
Equality holds for $n=2$. In particular, $\FF_{\infty}$ is exact on the category of finite length smooth representations of $\GL(2)$.
\end{cor}

\noindent {\it{Proof}}. Consider the Bernstein-Zelevinsky filtration of $V|_P$, as introduced at the very bottom of p. 48 in [BZ],
$$
V_{n-1}\subset V_{n-2}\subset \cdots \subset V_0=V
$$
The completely non-degenerate part $V^{nd}=V_{n-1}$ is a direct sum of $\dim V_{N,\psi}$ copies of the Kirillov representation
$\ind_N^P(\psi)$, and $V/V^{nd}$ is degenerate. As just shown,
$\FF_{\infty}(\ind_N^P(\psi))$ is one-dimensional, and therefore $\FF_{\infty}(V^{nd})$ has dimension $\dim V_{N,\psi}$. Since $\FF_{\infty}$ certainly is at least left-exact, the inequality follows.

\medskip

Now suppose $n=2$. The above filtration has length two, so $V/V^{nd}$ can be assumed to be given by a character $\chi$. Now, $\chi^{P(\OO)}=0$ unless $\chi$ is unramified, in which case $U^{(1)}$ acts on $\chi^{P(\OO)}=A$ by multiplication by $\chi(\varpi)\neq 0$. Ergo, $\FF_{\infty}(\chi)=0$, and therefore
$\FF_{\infty}(V)$ and $V_{N,\psi}$ have the same dimension. Since the functor $(-)_{N,\psi}$ is known to be exact (Proposition 2.35 in [BZ]), by a dimension-count we deduce that $\FF_{\infty}$ is exact too. Note that $V_{N,\psi}$ is finite-dimensional if $V$ has finite $G$-length -- see Corollary 5.22 on p. 51 in [BZ]. $\square$

\medskip

For $n>2$ there are more intricate non-generic representations $V$ of $\GL(n)$ than just characters, and it seems difficult to show that $\FF_{\infty}(V)=0$ for those. 

\medskip

As another application in the $\GL(2)$-case, and as a curiosum, we mention:

\begin{cor}
Let $V$ be an unramified principal series representation of $G=\GL_2(F)$ over $A=\bar{K}$ or $A=\bar{k}$, with Satake parameters $\{\alpha,\beta\}$, and consider the operator $U^{(1)}: V^{U_1(\varpi^r)}\rightarrow V^{U_1(\varpi^r)}$ for $r\geq 2$. Its Jordan normal form is
$$
U^{(1)}\sim J_{r-1}(0)\oplus J_1(q^{\frac{1}{2}}\alpha)\oplus J_1(q^{\frac{1}{2}}\beta),
$$
where $J_t(\lambda)$ is the $t \times t$ Jordan block with eigenvalue $\lambda$.
\end{cor}

\noindent {\it{Proof}}. From Proposition 4.4 in \cite{11}, we know that the generalized eigenspaces for $U^{(1)}$ on $V^{U_1(\varpi^r)}$ correspond to proper subsets
$\mathcal{T}\subset \{1,2\}$. The dimension is $\binom{r-1}{1-|\mathcal{T}|}$, which is $r-1$ or $1$ according to whether $\mathcal{T}$ is empty or not, and the pertaining eigenvalues are
$$
\lambda_{\varnothing}=0, \y \lambda_{\{1\}}=q^{\frac{1}{2}}\alpha, \y \lambda_{\{2\}}=q^{\frac{1}{2}}\beta.
$$
To see that $\lambda_{\varnothing}=0$ has only one Jordan block, we are to check the actual eigenspace is one-dimensional. But this is exactly $\mathcal{F}_r(V)$, which we know has dimension one by part (1) of Theorem 2. $\square$

\medskip

For instance, this shows that $U^{(1)}$ does {\it{not}} act semisimply on $V^{U_1(\varpi^r)}$ unless $r=2$ (a question left open in \cite{11}, p. 12). 
For $\GL(n)$ we have several operators $U^{(j)}$, and the {\it{simultaneous}} generalized eigespaces are indexed by subsets $\mathcal{T}$, as in the previous proof. This makes it difficult to pin down
the Jordan normal form of a specific operator $U^{(j)}$ -- when $n>2$.

\section{Notation -- the global setup}

Throughout, we will adopt the notation used in \cite{1}. Here we briefly recall a few essential definitions from their subsection 3.3, which will remain in force.

\subsection{Definite unitary groups}

Let $\ell>n>1$ be a prime, which splits in the imaginary quadratic field $E$. Fix a totally real field $F^+$, and introduce the CM-field $F=EF^+$. To be safe, we will always assume $F/F^+$ is unramified everywhere, although this may be removable. Choose a central division $F$-algebra $B$, of dimension $n^2$ over $F$, which satisfies the conditions:

\begin{itemize}
\item $B^{op} \simeq B \otimes_{E,c}E$,
\item $B_w\simeq M_n(F_w)$ at places $w$ above $w|_{F^+} \notin S(B)$,
\item $B_w$ is a \underline{division} algebra whenever $w|_{F^+} \in S(B)$.
\end{itemize}

\noindent Here $S(B)\neq \varnothing$ is a finite set of places of $F^+$, which split in $F$. We assume $S(B)$ is disjoint from $S_{\ell}$ (the places of $F^+$ above $\ell$), and in addition require that
$\#S(B)$ has the same parity as $\frac{n}{2}[F^+:\Q]$ when $n$ is even. We endow $B$ with an $F^+$-linear anti-involution $\star$ of the second kind ($\star|_F=c$), and
let $G=U(B,\star)_{/F^+}$ be the associated unitary group. We will always assume the pair $(B,\star)$ is chosen such that the following two bullets are fulfilled.

\begin{itemize}
\item $G(F_v^+)\simeq U(n)$ is {\it{compact}}, for all $v \in S_{\infty}$.
\item $G$ is {\it{quasi-split}} over $F_v^+$, for all $v \notin S(B)$.
\end{itemize}

\noindent By pinning down an order $\OO_B$, we may even define a model of $G$ over $\OO_{F^+}$, which will still be denoted by $G$. Suppose $v$ splits in $F$, say $v=ww^c$.
The choice of a place $w|v$ determines isomorphisms (unique up to conjugacy),
$$
\text{$i_w: G(F_v^+)\overset{\sim}{\longrightarrow} \GL_n(F_w)$, $\y$ $i_w: G(F_v^+)\overset{\sim}{\longrightarrow} B_w^{\times}$,} 
$$
according to whether $v$ lies outside $S(B)$ or not. We can (and will) even arrange that $i_w$ identifies $G(\OO_{F_v^+})$ with $\GL_n(\OO_{F_w})$ and $\OO_{B,w}^{\times}$ respectively.

\medskip

We will always choose our coefficient field $K$ large enough. That is, $K/\Q_{\ell}$ is a finite extension such that every embedding
$F \hookrightarrow \bar{K}$ actually maps into $K$. We let $\OO=\OO_K$ be its valuation ring, and $k=\OO/\lambda$ its residue field. Note,
$$
I_{\ell}=\Hom(F^+,K) \overset{\sim}{\longrightarrow} \bigsqcup_{v \in S_{\ell}} \Hom_{cts.}(F_v^+,K),
$$ 
and similarly for $F$. If for each $v \in S_{\ell}$ we have marked a place $\tilde{v}$ of $F$ above it, this collection $\{\tilde{v}\}$ defines a subset $\tilde{I}_{\ell}\subset \Hom(F,K)$ 
such that $\tilde{I}_{\ell} \overset{\sim}{\longrightarrow} I_{\ell}$ via the obvious restriction map. 

\subsection{Algebraic modular forms}

\subsubsection{Weights}

Let $\Z_+^n$ denote all $n$-tuples of integers $a=(a_1\geq \cdots \geq a_n)$. Each such $a$ is the highest weight of a unique irreducible algebraic representation
$\xi_a: \GL_n\rightarrow \GL(W_a)$ over $\Q$, which can even be defined over $\Z$. Let $M_a \subset W_a$ be a $\Z$-lattice. $\Wt_n$ is the subset of $(\Z_+^n)^{\Hom(F,K)}$ consisting of $a=(a_{\tau})_{\tau:F \rightarrow K}$ which satisfy the polarization condition, $a_{\tau c,i}=-a_{\tau,n+1-i}$. With each such $a$, we associate an irreducible representation,
$$
\xi_a: G(F^+\otimes_{\Q} \Q_{\ell}) \hookrightarrow \prod_{\tau \in \tilde{I}_{\ell}} \GL_n(K)\longrightarrow \GL(W_a)=\GL(\otimes_{\tau \in \tilde{I}_{\ell}}W_{a_{\tau}}).
$$
(Here we abuse notation, and view $W_{a_{\tau}}$ over $K$.) There is a natural $\OO$-lattice $M_a \subset W_a$, stable under the action of 
the group $G(\OO_{F^+}\otimes_{\Z}\Z_{\ell})$.

\subsubsection{Types}

For each place $v=ww^c \in S(B)$, suppose we are given an absolutely irreducible representation $\rho_v: G(F_v^+)\rightarrow \GL(M_{\rho_v})$, with open kernel, on a finite free $\OO$-module $M_{\rho_v}$. Via Jacquet-Langlands, $\rho_v \circ i_w^{-1}$ corresponds to a generalized Steinberg representation $\Sp_{n/m_v}(\pi_w)$ of $\GL_n(F_w)$; here $\pi_w$ is a supercuspidal representation of $\GL_{n/m_v}(F_w)$, which in turn corresponds to $\tilde{r}_w: \Gamma_{F_w}\rightarrow \GL_{n/m_v}(\bar{K})$ under local 
Langlands (suitably normalized). We will always assume $K$ is large enough, so that $\tilde{r}_w$ in fact maps into $\GL_{n/m_v}(\OO)$.

\bigskip

{\it{Remark}}. In \cite{1} they also operate with a set of places $R$, and for each $w|v$ with $v \in R$, specify characters $\chi_w:k(w)^{\times n}\rightarrow \OO^{\times}$, which play the role of eigen-characters for the action of Iwahori $\text{Iw}(w)$ on the $\text{Iw}_1(w)$-invariants. However, in their modularity lifting theorems (such as their Theorem 5.4.1) they eventually take $R=\varnothing$. For simplicity we will take $R$ to be empty throughout.

\bigskip

Choose a level-subgroup $U\subset G(\A_{F^+}^{\infty})$, which we assume factors as $\prod_v U_v$. We often assume $U$ is sufficiently small; by which we mean at least some 
$U_v$ contains no non-trivial elements of finite order. Moreover, we assume $U_v \subset G(\OO_{F_v^+})$ for every $v \in S_{\ell}$ -- so that $U_v$ acts on $M_a$. Introduce
$$
M_{a,\{\rho_v\},\varnothing}=M_a \otimes_{\OO} ({\bigotimes}_{v \in S(B)} M_{\rho_v}),
$$
an $\OO$-module with commuting actions of $G(\OO_{F^+}\otimes \Z_{\ell})$ and $G(F_{S(B)}^+)$. We include the subscript $\varnothing$ to emphasize we take 
$R$ empty; to stick with the notation in \cite{1}. Now, for any $\OO$-algebra $A$, let $S_{a,\{\rho_v\},\varnothing}(U,A)$ be the set of functions,
$$
\text{$f: G(F^+)\backslash G(\A_{F^+}^{\infty})\rightarrow M_{a,\{\rho_v\},\varnothing}\otimes_{\OO}A$, $\y$ $f(gu)=u_{S_{\ell}\cup S(B)}^{-1}f(g)$,}
$$
for all $u\in U$. When $U$ is sufficiently small this is a finite free $A$-module. We will make frequent use of its automorphic description over $\bar{K}\simeq \C$,
$$
S_{a,\{\rho_v\},\varnothing}(U,\OO)\otimes\bar{K}\simeq \bigoplus_{\pi} m_{\pi} \cdot (\otimes_{v \in S(B)} \Hom_{U_v}(\rho_v^{\vee},\pi_v)) \otimes (\otimes_{v \notin S(B)}\pi_v^{U_v}),
$$
where $\pi$ runs over automorphic representations of $G(\A_{F^+})$ with $\pi_{\infty}\simeq \xi_a^{\vee}$.

\subsection{Hecke algebras and Galois representations}

Let $T \supset S_{\ell}\cup S(B)$ be a finite set of split places of $F^+$. Suppose $U_v=G(\OO_{F_v^+})$ for all split places $v \notin T$. For each of the two divisors $w|v$, there are Hecke operators $T_w^{(j)}$ on modular forms, where $j=1,\ldots,n$, defined as the coset operators
$$
\text{$T_w^{(j)}=[\GL_n(\OO_{F_w})\zeta_{w,j}\GL_n(\OO_{F_w})]$, $\y$ $\zeta_{w,j}=\begin{pmatrix} \varpi_w\cdot I_{j} & \\ & I_{n-j} \end{pmatrix}$.}
$$
(See p. 99 of \cite{1} for more details; they act via the fixed isomorphism $i_w$.) We let
$$
\T_{{a,\{\rho_v\},\varnothing}}^T(U)=\OO[T_w^{(1)},\ldots,T_w^{(n)\pm 1}]_{w|v \notin T}\subset \End_{\OO}(S_{a,\{\rho_v\},\varnothing}(U,\OO))
$$
denote the $\OO$-subalgebra generated by the $T_w^{(j)}$ (along with $T_w^{(n)-1}$) where $w$ ranges over places of $F$ such that 
$v=w|_{F^+}\notin T$ splits in $F$. Thus $\T_{{a,\{\rho_v\},\varnothing}}^T(U)$ is a reduced commutative algebra, finite free over $\OO$.

\bigskip

For each maximal ideal $\m \subset \T_{{a,\{\rho_v\},\varnothing}}^T(U)$ (whose residue field is a finite extension of $k$) there is an associated continuous semisimple Galois representation
$$
\bar{r}_{\m}: \Gamma_F=\Gal(\bar{F}/F)\rightarrow \GL_n(\T_{{a,\{\rho_v\},\varnothing}}^T(U)/\m),
$$
as in 3.4.2 of \cite{1}. Suppose $v=ww^c \notin T$; then $\bar{r}_{\m}$ is unramified at $w$, and the characteristic polynomial of $\bar{r}_{\m}(\Frob_w)$ is the $\GL_n$-Hecke polynomial. We call $\m$ {\it{non-Eisenstein}} when $\bar{r}_{\m}$ is absolutely irreducible -- in which case it has a natural continuous lifting $r_{\m}$ to the localization
$\T_{{a,\{\rho_v\},\varnothing}}^T(U)_{\m}$, as described in 3.4.4 of \cite{1}. (Again, we may have to enlarge $K$ to ensure $\m$ has residue field exactly $k$.)

\section{Modularity lifting in the non-minimal case}

The proof of $R=\T$ mimics Wiles's proof for $\GL(2)$: It comes down to verifying a certain numerical criterion (comparing the size of the tangent space with that of a congruence module). In the minimal case, where no ramification is allowed, this can be done by a certain intricate patching argument -- the construction of Taylor-Wiles systems. The general case is then reduced to the minimal case by controlling level-raising congruences on the Hecke side -- which at the end comes down to Ihara's lemma.

\subsection{Ihara's lemma for unitary groups}

Throughout, we will admit the following conjecture. ("Ihara's lemma".)

\begin{cn}
Let $U \subset G(\A_{F^+}^{\infty})$ be a sufficiently small subgroup, and suppose
\begin{itemize}
\item $v \in T-(S_{\ell}\cup S(B))$ is a place where $U_v=G(\OO_{F_v^+})$,
\item $\m \subset \T_{a,\{\rho_x\},\varnothing}^T(U)$ is a non-Eisenstein maximal ideal,
\item $f \in S_{a,\{\rho_x\},\varnothing}(U,\bar{k})[\m]$ is an eigenform.
\end{itemize}
Then every irreducible $\bar{k}[G(F_v^+)]$-submodule $\pi$, 
$$
\pi \subset \la G(F_v^+)f \ra \subset S_{a,\{\rho_x\},\varnothing}(U^v,\bar{k}),
$$
is generic.
\end{cn}

In fact, one can deduce it from the special case of trivial data ($a=0$ and all $\rho_x=1$), which is Conjecture I on p. 145 of \cite{1}. Indeed, we may shrink $U$ at the places 
$x \in S_{\ell}\cup S(B)$, and assume $U_x$ is pro-$\ell$. It therefore acts unipotently on $M_a \otimes \bar{k}$. Thus, up to semisimplification,  
$S_{a,\{\rho_x\},\varnothing}(U^v,\bar{k})$ is just a direct sum of finitely many copies of $S_{0,\{1\},\varnothing}(U^v,\bar{k})$. Compare with Lemma 5.3.2 in \cite{1}.

\bigskip

In \cite{1} they propose a {\it{stronger}} conjecture, which they do not need however. The key difference with Conjecture 1 is that $\pi$ is not assumed to sit in a cyclic submodule generated by a $G(\OO_{F_v^+})$-spherical eigenform $f$. We restate it below.

\begin{cn}
Let $U \subset G(\A_{F^+}^{\infty})$ be a sufficiently small subgroup, and suppose
\begin{itemize}
\item $v \in T-(S_{\ell}\cup S(B))$ is a place (split in $F$, but no restriction on $U_v$),
\item $\m \subset \T_{a,\{\rho_x\},\varnothing}^T(U)$ is a non-Eisenstein maximal ideal.
\end{itemize}
Then every irreducible $\bar{k}[G(F_v^+)]$-submodule $\pi \subset S_{a,\{\rho_x\},\varnothing}(U^v,\bar{k})_{\m}$ is generic.
\end{cn}

\noindent This is Conjecture II on p. 146 in \cite{1} (see also their Conjecture B in the introduction of loc. cit.) -- where they take trivial data for simplicity, and without loss of generality (see their Lemma 5.3.2, again). 

\bigskip

\noindent At the very end of this paper (section 5.5 below) we will appeal to Conjecture 3 at a point where Conjecture 2 is insufficient. It turns out a slightly weaker statement is enough for our purposes; namely that any $\pi\subset S_{a,\{\rho_x\},\varnothing}(U^v,\bar{k})_{\m}[\m]$ is generic. We will apply this to a closely related group $G'=U(B^{op},\star)$.

\begin{lem}
The conjectures above hold for $n=2$.
\end{lem}

\noindent {\it{Proof}}. This is basically just the strong approximation theorem for $G^{der}$. If $\pi$ is non-generic, it is one-dimensional, $\pi=\bar{k}\cdot f$ where $f$ 
factors through the torus $\ker(N_{F/F^+})$. One easily deduces $\m$ must be Eisenstein. See 5.3.1 in \cite{1}. $\square$

\bigskip

Conjectures 2 and 3 are still open for $n>2$. The $R=\T$ theorems in \cite{1} rely on 2 (in the non-minimal case). Taylor has since adapted a variant of the patching argument, due to Kisin,
which avoids Ihara's lemma -- but proves the a priori weaker result $R^{red}=\T$. As is discussed in the introduction of \cite{1}, the stronger result $R=\T$ would still be of interest,
and have potential applications to special value formulas for adjoint $L$-functions. One of the goals of this paper is to recast the "Ihara" conjectures 2 and 3 in a more modern and robust framework; that of essentially AIG representations, which play a fundamental role in "local Langlands in families" -- due to Emerton and Helm.

\subsection{$R=\T$ in the non-minimal case}

Here we follow the discussion in section 5.4 of \cite{1} rather closely. Thus, additionally, we assume $F/F^+$ is unramified everywhere, and that $\ell$ is unramified in $F^+$.

\bigskip

We choose a finite set of auxiliary primes $S_a\neq \varnothing$, which all split in $F$. We assume $S_a$ is disjoint from $S_{\ell}\cup S(B)$. Moreover, if $v\in S_a$ lies above the rational prime $p$, then $[F(\zeta_p):F]>n$. The sole purpose of $S_a$ is to ensure that the levels we work with below are sufficiently small. 

\bigskip

Finally, $S$ is a finite set of split places, disjoint from $S_{\ell}\cup S(B) \cup S_a$, consisting of banal places $v$. That is, $\ell \nmid \#\GL_n(k(v))$ for all $v \in S$. All these places taken together gives us
$$
T=S \cup S_{\ell}\cup S(B) \cup S_a.
$$
For each $v \in T$, we once and for all choose a place $\tilde{v}$ of $F$ dividing it, and let $\tilde{T}$ be the collection of all these.

\bigskip

Now, for each subset $\Sigma \subset S$, we introduce a level-subgroup $U(\Sigma)$. (Here our notation differs slightly from \cite{1}, in that we use $\Sigma$ instead of the somewhat cumbersome $S_1$.) It is defined component-wise, by
$$
U(\Sigma)_v=
\begin{cases}
\text{hyperspecial}, & \text{$v$ non-split},\\
G(\OO_{F_v^+}), & \text{$v \notin \Sigma \cup S_a$ split},\\
i_{\tilde{v}}^{-1}U_1(\tilde{v}^n), & v \in \Sigma, \\
i_{\tilde{v}}^{-1} U(\tilde{v}), & v \in S_a.
\end{cases}
$$
Here $U(\tilde{v})$ is the kernel of the reduction map to $\GL_n(\OO_{F_{\tilde{v}}}/\varpi_{\tilde{v}})$, and $U_1(\tilde{v}^n)$ consists of matrices whose last row is
$(0^{n-1},1)$ mod $\tilde{v}^n$. (The exponent $n$ is the {\it{same}} as in $\GL_n$.) Let $U=U(\varnothing)$ -- which is maximal away from $S_a$.

\bigskip

We assume that the weight $a=(a_{\tau})_{\tau:F \rightarrow K} \in \Wt_n$ is in the "Fontaine-Laffaille" range. That is,
satisfies the smallness-condition:
$$
\ell-n>a_{\tau,1}\geq \cdots \geq a_{\tau,n}\geq 0
$$  
for all $\tau \in \tilde{I}_{\ell}$. (That is, those $\tau$ which define a place in $\tilde{S}_{\ell}$.) Furthermore, there is a small technical hypothesis on the types $\{\rho_v\}_{v \in S(B)}$. Namely, $\tilde{r}_{\tilde{v}}\otimes_{\OO}k$ should be absolutely irreducible, and not isomorphic to any cyclotomic twist $(\tilde{r}_{\tilde{v}}\otimes_{\OO}k)(i)$ for
$i=1,\ldots,m_v$.

\bigskip

Now, let $\m \subset \T_{a,\{\rho_x\},\varnothing}^T(U)$ be a non-Eisenstein maximal ideal, with residue field $k$, and let $\bar{r}_{\m}$ be its associated Galois representation.
We will assume throughout that the "modularity lifting criteria" are fulfilled: 

\begin{itemize}
\item $\bar{r}_{\m}(\Gamma_{F^+(\zeta_{\ell})})$ is big (in the sense of 2.5.1 on p. 55 of \cite{1}). 
\item $\bar{r}_{\m}$ is unramified at every $v \in S_a$, and $(\text{ad} \bar{r}_{\m})(1)^{\Gamma_{F_{\tilde{v}}}}=0$.
\end{itemize}

\noindent For $\Sigma \subset S$, consider the deformation problem $\mathcal{S}_{\Sigma}$ as in \cite{1}, where we only point out that one takes crystalline lifts at $v \in S_{\ell}$, and unramified lifts at $v \in S-\Sigma$. Let $R_{\Sigma}^{univ}$ denote the corresponding universal deformation ring of $\bar{r}_{\m}$.

\bigskip

Pull back $\m$ to a maximal ideal of $\T_{a,\{\rho_x\},\varnothing}^T(U(\Sigma))$, and consider the localized module,
$$
X_{\Sigma}=S_{a,\{\rho_x\},\varnothing}(U(\Sigma),\OO)_{\m,\n},
$$
where 
$$
\n=(\lambda,U_{\tilde{v}}^{(1)},\ldots,U_{\tilde{v}}^{(n-1)})_{v \in \Sigma}.
$$
The image of $\T_{a,\{\rho_x\},\varnothing}^T(U(\Sigma))_{\m}$ in $\End(X_{\Sigma})$ is then a complete local Noetherian $\OO$-algebra $\T_{\Sigma}$, which is reduced, and finite free over $\OO$. By composition, $r_{\m}$ gives a natural lifting of $\bar{r}_{\m}$ to $\T_{\Sigma}$, of type $\mathcal{S}_{\Sigma}$, and consequently there is a natural surjective map,
$$
R_{\Sigma}^{univ} \longrightarrow \T_{\Sigma}.
$$
One of the main results of \cite{1} shows this is an isomorphism if we admit the conjectural analogue of Ihara's lemma recalled previously.

\begin{thm}
Assume Conjecture 1 is true. Then $R_{\Sigma}^{univ} \simeq \T_{\Sigma}$. In addition, $X_{\Sigma}$ is a free module over $\T_{\Sigma}$. (Both hold unconditionally if $\Sigma=\varnothing$ or $n=2$.)
\end{thm}

\noindent {\it{Proof}}. This is Theorem 5.4.1 on p. 155 in \cite{1}, except that the freeness is not mentioned explicitly. However, it follows directly from the proof combined with Diamond's numerical freeness criterion (Theorem 2.4 in \cite{2}). $\square$

\bigskip

(We remark that the freeness of $X_{\varnothing}$ is mentioned in Theorem 3.5.1 of \cite{1}, and the freeness of $X_{\Sigma}$ is one of the hypotheses of their Lemma 5.3.3.) 

\bigskip

The overall strategy of the proof follows that of Wiles. Firstly one proves it in the minimal case, $\Sigma=\varnothing$, by patching together modules and algebras of auxiliary level 
(the Taylor-Wiles method). Secondly, one shows that $R_{\Sigma}^{univ}$ and $\T_{\Sigma}$ have the same growth rate as one enlarges $\Sigma$. To check this for the Hecke algebra, one needs to control the size of a congruence module, measuring level-raising congruences -- this is where Conjecture 1 enters.

\subsection{The rank of $X_{\Sigma}$ is one}

When $S_a=\varnothing$, one can make the exact same definitions as in the previous section, and get a $\T_{\Sigma}$-module $X_{\Sigma}$. The problem becomes that $U=U(\varnothing)$ is {\it{maximal}} everywhere, and hence not sufficiently small. However, freeness still holds.

\begin{lem}
Admitting Conjecture 1, $X_{\Sigma}$ is free over $\T_{\Sigma}$; even when $S_a=\varnothing$.
\end{lem}

\noindent {\it{Proof}}. Pick a split place $u$, above a large enough rational prime, which is banal for $\ell$. That is, $\ell \nmid \#\GL_n(k(u))$. Let $S_a=\{u\}$, and define corresponding
$X_{\Sigma}'$ and $\T_{\Sigma}'$ as in the previous section, of level $U'(\Sigma)$ contained in $U(\Sigma)$ -- they differ only at the place $u$, where $U'(\Sigma)_u\simeq 
U(\tilde{u})$, whereas $U(\Sigma)_u=G(\OO_{F_u^+})$. By $R=\T$, we know that $X_{\Sigma}'$ is free over $\T_{\Sigma}'$. Now, $X_{\Sigma} \hookrightarrow X_{\Sigma}'$, and we view both as $\T_{\Sigma}'$-modules via the restriction map $\T_{\Sigma}' \twoheadrightarrow \T_{\Sigma}$. The usual averaging-idempotent $e_{U(\tilde{u})}$ is defined over
$\OO$, since $u$ is banal for $\ell$, and obviously commutes with the Hecke-action away from $T$. Consequently, $X_{\Sigma} \hookrightarrow X_{\Sigma}'$ admits a retraction, showing that $X_{\Sigma}$ sits as a direct summand of $X_{\Sigma}'$, which is free. In other words, $X_{\Sigma}$ is at least a projective $\T_{\Sigma}'$-module. Since $\T_{\Sigma}'$ is a local ring, a result of Kaplansky tells us $X_{\Sigma}$ is in fact free over $\T_{\Sigma}'$. (In particular, $\T_{\Sigma}' \simeq \T_{\Sigma}$.) $\square$

\bigskip

From now on, we will assume $S_a=\varnothing$. Moreover, we will make the following multiplicity one hypothesis for the (non quasi-split) unitary group $G$.

\begin{hyp}
Every automorphic representation $\pi$ of $G(\A_{F^+})$ occurs with multiplicity $m_{\pi}=1$.
\end{hyp}

We are optimistic that the trace formula experts have checked this, but have been unsuccessful in finding a precise reference. At the very least we need this for stable $\pi$ -- by which we mean its Galois representation $r_{\pi}$ is irreducible.

\medskip

\noindent {\it{Remark}}. For quasi-split $U(n)$, multiplicity one was proved for $n=2$ and $n=3$ by Rogawski. See Theorem 11.5.1(c) and Theorem 13.3.1 in \cite{12}. (For $n=2$, see also Corollary 5 on p. 726 of \cite{13}.) For general $n$, multiplicity one presumably follows from the much more general results of Mok (which mimics work of Arthur for symplectic and orthogonal groups) -- again, in the quasi-split case. See \cite{14}. For general unitary groups $G$ (that is, inner forms of $U(n)$, possibly non-quasi-split), $m_{\pi}=1$ should follow from work in progress of Kaletha, Shin, and White. In Chapter 14 of \cite{12}, Rogawski deals with general $G$, in two or three variables, and establishes the analogue of the Jacquet-Langlands correspondence -- from which one deduces multiplicity one. To summarize, at least the multiplicity one hypothesis is fulfilled for $n=2$ and $n=3$.

\begin{prop}
Admitting Hypothesis 1, $X_{\Sigma}$ is free of rank \underline{one} over $\T_{\Sigma}$ (assuming $S_a=\varnothing$).
\end{prop}

\noindent {\it{Proof}}. Let $r$ be the rank of $X_{\Sigma}$. We compute it by tensoring with $\bar{K}$. Note,
$$
X_{\Sigma}\otimes_{\OO}\bar{K}\simeq  
$$
$$
\bigoplus_{\pi} (\otimes_{v \in S(B)} \Hom_{U_v}(\rho_v^{\vee},\pi_v)) \otimes (\otimes_{v \in S_{\ell}\cup S\backslash\Sigma}\pi_v^{U_v}) \otimes (\otimes_{v\notin T}\pi_v^{U_v})_{\m} \otimes 
(\otimes_{v \in \Sigma}\pi_v^{U(\Sigma)_v})_{\n},
$$
where $\pi$ runs over automorphic representations of $G(\A_{F^+})$ with $\pi_{\infty}\simeq \xi_a^{\vee}$. We have omitted $m_{\pi}$ in the decomposition, making use of Hypothesis 1. This is a free module of rank $r$ over
$$
\T_{\Sigma}\otimes_{\OO}\bar{K} \simeq \bar{K}^{\Hom_{\OO-alg.}(\T_{\Sigma},\bar{K})}.
$$
The tensor product $\otimes_{v\in S(B)}$ is one-dimensional over $\bar{K}$. Indeed $U_v \simeq \OO_{B,\tilde{v}}^{\times}$, so $G(F_v^+)$ is generated by $U_v$ and the center.
Thus $\pi_v|_{U_v}$ remains irreducible, and hence isomorphic to $\rho_v^{\vee}$. Now use Schur. Secondly, the next factor $\otimes_{v \in S_{\ell}\cup S\backslash\Sigma}$ is one-dimensional since $U_v\simeq \GL_n(\OO_{F_{\tilde{v}}})$. For the same reason, $\otimes_{v \notin T}$ is one-dimensional -- $U_v$ being hyperspecial for $v \notin T$.
Finally, the $\n$-localization of $\otimes_{v\in \Sigma}$ is one-dimensional -- and this is the crux of the mater. We give the proof in separate paragraphs:

\bigskip

For each $\pi$ contributing to $X_{\Sigma}\otimes_{\OO}\bar{K}$, we must show $\LL_n(\pi_v)\simeq \bar{K}$ for all $v \in \Sigma$. Here we think of $\pi_v$ as a representation of
$\GL_n(F_{\tilde{v}})$ via our choice $i_{\tilde{v}}$. By 3.3.4 in \cite{1}, we know how to attach a Galois representation $r_{\pi}: \Gamma_F \rightarrow \GL_n(\OO')$, where $\OO'\supset \OO$, which is now known to satisfy local-global compatibility at all finite primes (although we only need the result away from $\ell$). Let $r_{\pi,\tilde{v}}=r_{\pi}|_{\Gamma_{F_{\tilde{v}}}}$ be shorthand notation. Via local Langlands, $r_{\pi,\tilde{v}} \leftrightarrow \pi_v$, and this is a {\it{generic}} representation: Indeed, $\bar{r}_{\pi}\simeq \bar{r}_{\m}$, and the latter is irreducible (since $\m$ is assumed to be non-Eisenstein). A fortiori, $r_{\pi}$ must be irreducible, which means the base change $BC_{F/F^+}(\pi)$
is cuspidal, and therefore (globally) generic. In particular, $\pi_v$ is generic. Thus local Langlands coincides with Breuil-Schneider's "generic" local Langlands in this case. 

\bigskip

Now, as a key input, we use the (modified) mod $p$ local Langlands correspondence of Emerton and Helm, cf. Theorem 5.1.5 in \cite{5} -- which builds on the semisimple correspondence of Vigneras. Since $r_{\pi,\tilde{v}}$ is a lift of $\bar{r}_{\m,\tilde{v}}\otimes_k k'$, there is a unique (up to homothety) $\GL_n(F_{\tilde{v}})$-invariant $\OO'$-lattice
$L \subset \pi_v$ such that $L\otimes k'$ is essentially AIG, and embeddable in $\bar{\pi}(\bar{r}_{\m,\tilde{v}})\otimes_k k'$. This is basically just paraphrasing part (2) of 
5.1.5 in \cite{5}. Since $\LL_n$ commutes with extensions of scalars, it suffices to show $\LL_n(L)\simeq \OO'$. Or in turn, that $\LL_n(L\otimes k')\simeq k'$. To do this, recall 
part (8) of 5.1.5 in \cite{5}, which says that
all constituents of $\bar{\pi}(\bar{r}_{\m,\tilde{v}})$ (and hence of $L \otimes k'$) have the same supercuspidal support. Namely,
$$
\text{$|\cdot|^{\frac{n-1}{2}}\cdot \{\bar{\chi}_1,\ldots,\bar{\chi}_n\}$, $\y$ $\bar{r}_{\m,\tilde{v}}^{ss}\simeq \begin{pmatrix}\bar{\chi}_1 & & \\ & \ddots & \\ & & \bar{\chi}_n\end{pmatrix}$.}
$$
We remind ourselves that $\m$ comes from full level $U=U(\varnothing)$, so $\bar{r}_{\m}$ is unramified away from $\ell$; in particular at $\tilde{v}$. Thus the $\bar{\chi}_i$ are 
unramified characters $\Gamma_{F_{\tilde{v}}}\rightarrow \bar{k}^{\times}$, which we tacitly view as unramified characters of $F_{\tilde{v}}^{\times}$ via the local Artin map. 
What is important to us, is that every constituent of $L \otimes k'$ is a constituent of an unramified principal series $i_B(\bar{\chi_i})$ (ignoring the twist for the sake of clarity). 
We now invoke our local results. We showed that only the generic constituent of $i_B(\bar{\chi_i})$ has a nonzero $\LL_n$, which was shown to be one-dimensional. Therefore, since $L\otimes k'$ has a unique generic constituent (its socle), being essentially AIG, we deduce from the exactness of $\LL_n$ that $\LL_n(L\otimes k')$ is one-dimensional, as we wished. Consequently, so is $\LL_n(\pi_v)$.

\bigskip

It remains to show that a contributing $\pi$ is completely determined by its Hecke eigensystem $\phi_{\pi}: \T_{\Sigma} \rightarrow \bar{K}$. Certainly $\pi_v$ is determined for all $v \notin T$ split in $F$. Let $\Pi=BC_{F/F^+}(\pi)$ be the base change to $\GL_n(\A_F)$; whose existence is guaranteed by Proposition 3.3.2 of \cite{1}, for instance -- which quotes results of Clozel and Labesse. Then $\Pi_w$ is determined at all places $w$ such that $w|_{F^+}\notin T$ splits. Since $\Pi$ is isobaric, the main result (Theorem A) of \cite{3} shows $\Pi$ is uniquely determined: If $\Pi$ and $\Pi'$ are isobaric automorphic representations of $\GL_n(\A_F)$, such that $\Pi_w \simeq \Pi_w'$ for all but finitely many degree-one places $w$ of $F$ (over $F^+$) then $\Pi\simeq \Pi'$ -- alternatively one could work with the Galois representation $r_{\pi}$ from 3.3.4 in \cite{1}, and just use 
Tchebotarev. We still have to argue that $\pi$ is determined by its base change $\Pi$. The potential problem is what happens at the non-split places $v$ of $F^+$, where we must show the $L$-packets are singletons. However, at non-split $v$ we take $U_v$ to be a hyperspecial subgroup -- so $\pi_v$ is $U_v$-unramified. For unramified representations, base change is injective: We refer to Corollary 4.2 on p. 17 of \cite{4}. We conclude $r=1$. $\square$

\section{End of the proof of the main result}

\subsection{More notation and finite generation}

We keep the notation already introduced. For $U^{\Sigma}=\prod_{v\notin \Sigma}U_v$, we consider
$$
S_{a,\{\rho_v\},\varnothing}(U^{\Sigma},\OO)=\underset{U_{\Sigma}}{\varinjlim} S_{a,\{\rho_v\},\varnothing}(U_{\Sigma}U^{\Sigma},\OO),
$$
consisting of smooth $M_{a,\{\rho_v\},\varnothing}$-valued functions $f$ on $G(\A_{F^+}^{\infty})$, with the usual invariance properties. It carries a faithful action of the Hecke algebra,
$$
\T_{a,\{\rho_v\},\varnothing}^T(U^{\Sigma})=\underset{U_{\Sigma}}{\varprojlim}\T_{a,\{\rho_v\},\varnothing}^T(U_{\Sigma}U^{\Sigma}),
$$
and we pull back $\m$ to a maximal ideal in here (non-Eisenstein, and with residue field $k$). Then we introduce the module
$$
H=H_{\Sigma}=S_{a,\{\rho_v\},\varnothing}(U^{\Sigma},\OO)_{\m}.
$$
Thus $H \otimes_{\OO}\bar{K}$ is an admissible representation of $G(F_{\Sigma}^+)$, which we always identify with $\prod_{v\in \Sigma}\GL_n(F_{\tilde{v}})$. The identification obviously depends on the choices $\tilde{\Sigma}$ and $i_{\tilde{v}}$. As a representation of $G(F_{\Sigma}^+)$, again admitting Hypothesis 1 ($m_{\pi}=1$), 
$$
H \otimes_{\OO}\bar{K} \simeq \bigoplus_{\pi} \pi_{\Sigma},
$$
where $\pi$ runs over all automorphic representations of $G(\A_{F^+}^{\infty})$, of infinity type $\pi_{\infty}\simeq \xi_{a}^{\vee}$, types $\rho_v^{\vee}$ at $v \in S(B)$, with non-trivial $U_v$-invariants for $v \in S_{\ell}\cup S\backslash \Sigma$, $U_v$-unramified for $v \notin T$, and with $\bar{r}_{\pi}\simeq \bar{r}_{\m}$. Here $\pi_{\Sigma}=\otimes_{v\in \Sigma}\pi_v$.

\begin{lem}
$H \otimes_{\OO}\bar{K}$ has finite length as a $G(F_{\Sigma}^+)$-representation.
\end{lem}

\noindent {\it{Proof}}. Equivalently, by Howe's Theorem (see Theorem 4.1 on p. 37 in \cite{6}), we are asking that $H \otimes_{\OO}\bar{K}$
is finitely generated. Since $\pi_{\Sigma}$ is irreducible (hence cyclic), we are to show there are only finitely many $\pi$ contributing to 
$H \otimes_{\OO}\bar{K}$ in the above decomposition: Any such $\pi$ has a fixed weight $a$, and $\pi_v^{U_v}\neq 0$ for 
$v \in S_{\ell}\cup S\backslash \Sigma$, and for $v \notin T$. Moreover, $\pi_v^{\ker(\rho_v)}\neq 0$ for $v \in S(B)$. It remains to check that the condition
$\bar{r}_{\pi}\simeq \bar{r}_{\m}$ forces the conductor of $\pi_v$ to be bounded for $v \in \Sigma$. This follows from an observation of Livne (and independently Carayol). Indeed, 
since $r_{\pi}$ is a lift of $\bar{r}_{\m}$, it follows from Proposition 1.1 on p. 135 of \cite{7} that $r_{\pi,\tilde{v}}^{ss}$ and $\bar{r}_{\pi,\tilde{v}}$ have the same Swan condutors
($\Sigma \cap S_{\ell}=\varnothing$), for $v \in \Sigma$. As a result, the Artin conductor of $r_{\pi,\tilde{v}}$ is at most that of $\bar{r}_{\pi,\tilde{v}}$ plus $n$. We conclude, from preservation of $\epsilon$-factors, that the conductor of $\pi_v$ is at most $n$ -- since $\bar{r}_{\m,\tilde{v}}$ is unramified. $\square$

\subsection{Lattices and strong Brauer-Nesbitt}

Let $\pi$ be an automorphic representation appearing in the decomposition of $H \otimes_{\OO}\bar{K}$. Fix a place $v \in \Sigma$, for now. Since $r_{\pi,\tilde{v}}$ is a lift of
$\bar{r}_{\m,\tilde{v}}\otimes_k \bar{k}$, there is an $\OO_{\bar{K}}$-lattice $L_{\pi_v}\subset \pi_v$ such that $L_{\pi_v}\otimes \bar{k}$ is essentially AIG, and embeddable in 
$\bar{\pi}(\bar{r}_{\m,\tilde{v}})\otimes_k \bar{k}$. Moreover, such an $L_{\pi_v}$ is unique up to homothety. This is part of Theorem 5.1.5 in \cite{5}, and was already used previously in the proof of Proposition 1 above. (The analogue holds over a finite extension $k' \supset k$, the residue field of a discrete valuation ring $\OO' \supset \OO$, and 
we may sometimes think of $L_{\pi_v}$ as being defined over $\OO'$ -- enlarged if necessary.) Let $L_{\pi_{\Sigma}}=\otimes_{v \in \Sigma}L_{\pi_v}$.

\begin{lem}
$L=L_{\Sigma}=\oplus_{\pi} L_{\pi_{\Sigma}}$ is a finitely generated $G(F_{\Sigma}^+)$-stable $\OO_{\bar{K}}$-lattice in $H \otimes_{\OO}\bar{K}$. Furthermore, 
its reduction $L\otimes \bar{k}$ is of finite length, and $(L\otimes \bar{k})^{ss}$ is independent of the choice of (any) lattice $L$. In particular,
$$
(L\otimes \bar{k})^{ss} \simeq (H \otimes \bar{k})^{ss}.
$$
\end{lem}

\noindent {\it{Proof}}. This follows from the strong Brauer-Nesbitt principle of Vigneras (which is Theorem 1 in \cite{8}). It applies since we know by the previous lemma that 
$H \otimes_{\OO}\bar{K}$ is of finite length. $\square$

\bigskip

Note that the space on the right,
$$
H \otimes \bar{k}=S_{a,\{\rho_v\},\varnothing}(U^{\Sigma},\bar{k})_{\m},
$$
is a space of modular forms mod $\ell$. We infer that any $G(F_{\Sigma}^+)$-constituent hereof occurs in $L\otimes \bar{k}$ (with the same multiplicity), and therefore in
$\otimes_{v\in \Sigma}(L_{\pi_v}\otimes \bar{k})$, for some $\pi$. The latter tensor product can be embedded in $\otimes_{v\in \Sigma}\bar{\pi}(\bar{r}_{\m,\tilde{v}})$, after extending scalars to $\bar{k}$. As we have used before (in the proof of Proposition 1), all constituents of $\bar{\pi}(\bar{r}_{\m,\tilde{v}})$ appear in the same unramified principal series
$i_B(\bar{\chi}_{i,\tilde{v}})$. See 5.1.5 in \cite{5}. The upshot is that any constituent of $H \otimes \bar{k}$ occurs in (a tensor product of) unramified principal series
$\otimes_{v \in \Sigma} i_B(\bar{\chi}_{i,\tilde{v}})$.

\bigskip

Below it will be useful to have alternative ways to think of $\T_{\Sigma}$. For instance, we will often think of $H$ as a $\T_{\Sigma}$-module, via the following isomorphisms.

\begin{lem}
$\T_{a,\{\rho_v\},\varnothing}^T(U^{\Sigma})_{\m} \overset{\sim}{\longrightarrow} \T_{a,\{\rho_v\},\varnothing}^T(U(\Sigma))_{\m} \overset{\sim}{\longrightarrow} \T_{\Sigma}$.
\end{lem}

\noindent {\it{Proof}}. Suppose $t \in \T_{a,\{\rho_v\},\varnothing}^T(U^{\Sigma})_{\m}$ acts trivially on $X_{\Sigma}$. We must show $t$ acts trivially on $H$, so that $t=0$. It is enough to show $t$ acts trivially on $H \otimes_{\OO} \bar{K}$, which has an automorphic description, given before Lemma 3. Given a $\pi$ contributing to $H \otimes_{\OO} \bar{K}$, with eigensystem $\phi_{\pi}: \T_{a,\{\rho_v\},\varnothing}^T(U^{\Sigma}) \rightarrow \bar{K}$, we must show $\phi_{\pi}(t)=0$. This follows from our hypothesis on $t$ if we can verify $\pi$ contributes to $X_{\Sigma}\otimes_{\OO}\bar{K}$, as given in the beginning of the proof of Proposition 1. In other words, that $\LL_n(\pi_v)\neq 0$ for all $v \in \Sigma$.
It suffices to show $\LL_n(L_{\pi_v})\neq 0$. Or, in turn, that $\LL_n(L_{\pi_v}\otimes \bar{k})\neq 0$ -- all because $\LL_n$ commutes with extension of scalars. However, 
as noted above, $L_{\pi_v}\otimes \bar{k}$ is essentially AIG, and all its constituents appear in the same unramified principal series $i_B(\bar{\chi}_{i,\tilde{v}})$. In particular, its generic socle has a non-trivial $\LL_n$, as shown in our section on local preliminaries. $\square$

\subsection{Generic constituents of $H/\m H$}

As a first approximation to our main result, we have:

\begin{lem}
$H/\m H$ has a unique generic $G(F_{\Sigma}^+)$-constituent (occurring with multiplicity one). 
\end{lem}

\noindent {\it{Proof}}. First observe that $X_{\Sigma}=\LL_{n,\Sigma}(H)$, where by $\LL_{n,\Sigma}$ we mean the composition of all the functors $\LL_{n,\tilde{v}}$ for $v \in \Sigma$. It commutes with extension of scalars, so (viewing $H$ as a $\T_{\Sigma}$-module)
$$
\LL_{n,\Sigma}(H/\m H)=\LL_{n,\Sigma}(H \otimes_{\T_{\Sigma}} k)=\LL_{n,\Sigma}(H) \otimes_{\T_{\Sigma}} k=X_{\Sigma}\otimes_{\T_{\Sigma}} k =k,
$$
where, in the last step, we have used that $X_{\Sigma}$ is of rank one over $\T_{\Sigma}$. Since $\lambda \subset \m$, there is a surjection
$H/\lambda H \twoheadrightarrow H/\m H$, from which we deduce that every constituent of $H/\m H$ also occurs in $H/\lambda H=H \otimes_{\OO} k$. From Lemma 4, and the pertaining discussion, we see that such constituents appear in $\otimes_{v\in \Sigma} i_B(\bar{\chi}_{i,\tilde{v}})$. The unique constituent of $H/\m H$ with a nonzero $\LL_{n,\Sigma}$ must therefore be the unique generic constituent of this unramified principal series. $\square$

\subsection{Duality}

In Section 5.2 of \cite{1}, they define a Petersson pairing between automorphic forms on $G$ and those on a related group $G'=U(B^{op},\star)_{/F^+}$, isomorphic to $G$ via the inversion map $I: G \overset{\sim}{\longrightarrow} G'$. The choice of an order $\OO_B$ even yields an integral model of $G'$ over $\OO_{F^+}$, and $I$ is defined over $\OO_{F^+}$ as well. We choose isomorphisms $i_w^t$ as in the first few paragraphs on p. 143 in \cite{1}.

\bigskip

For each $a \in \Wt_n$, there is a natural pairing $\la,\ra_{a_{\tau}}$ on $W_{a_{\tau}}$, reflecting the fact that $\xi_{a_{\tau}}^{\vee}$ can be realized as $\xi_{a_{\tau}} \circ {^t(\cdot)}^{-1}$. We let $M_{a_{\tau}}'\subset W_{a_{\tau}}$ denote the $\Z$-lattice dual to $M_{a_{\tau}}$. Correspondingly, there is an irreducible algebraic representation,
$$
\xi_a':G'(F^+ \otimes_{\Q}\Q_{\ell}) \hookrightarrow \prod_{\tau \in \tilde{I}_{\ell}}\GL_n(K)\longrightarrow \GL(W_a)=\GL(\otimes_{\tau \in \tilde{I}_{\ell}}W_{a_{\tau}}),
$$
where we abuse notation and write $W_{a_{\tau}}$ instead of $W_{a_{\tau}}\otimes_{\Q}K$. By definition, the action of $G'(\OO_{F^+}\otimes_{\Z}\Z_{\ell})$ stabilizes the 
$\OO$-lattice $M_a'=\otimes_{\tau\in \tilde{I}_{\ell}}M_{a_{\tau}}'$. Furthermore, there is a perfect pairing 
$\la,\ra_a: M_a \times M_a' \rightarrow \OO$; contravariant in the sense that
$$
\la \xi_a(g)m,m'\ra_a=\la m, \xi_a'(I(g)^{-1})m'\ra_a
$$
for $g \in G(\OO_{F^+}\otimes_{\Z}\Z_{\ell})$. In particular, we identify $M_a' \simeq \Hom_{\OO}(M_a,\OO)$.

\bigskip

For each place $v=ww^c \in S(B)$, we introduce the finite free $\OO$-module $M_{\rho_v}'=\Hom_{\OO}(M_{\rho_v},\OO)$, equipped with a natural $G'(F_v^+)$-action 
$\rho_v'$. Namely,
$$
\la \rho_v(g)m,m'\ra=\la m, \rho_v'(I(g)^{-1}m')\ra
$$
for $g \in G(F_v^+)$. We consolidate all these data in the following $\OO$-module,
$$
M_{a,\{\rho_v\},\varnothing}'=M_a' \otimes_{\OO}({\bigotimes}_{v\in S(B)}M_{\rho_v}'),
$$
which carries commuting actions of $G'(\OO_{F^+}\otimes \Z_{\ell})$ and $G'(F_{S(B)}^+)$. Then for each compact open subgroup $U' \subset G'(\A_{F^+}^{\infty})$, 
and any $\OO$-algebra $A$, we define $S_{a,\{\rho_v\},\varnothing}'(U',A)$ to be the set of functions
$$
\text{$f': G'(F^+)\backslash G'(\A_{F^+}^{\infty})\rightarrow M_{a,\{\rho_v\},\varnothing}'\otimes_{\OO}A$, $\y$ $f'(gv)=v_{S_{\ell}\cup S(B)}^{-1}f'(g)$.}
$$
When $U$ is sufficiently small, and $U'=I(U)$, this is in duality with $S_{a,\{\rho_v\},\varnothing}(U,A)$
via the pairing defined at the bottom of p. 144 in \cite{1} (taking $\eta=1$),
$$
\la f,f'\ra={\sum}_{g \in G(F^+)\backslash G(\A_{F^+}^{\infty})/U} \la f(g),f'(I(g))\ra.
$$
(We keep assuming $U_v \subset G(\OO_{F_v^+})$ for $v \in S_{\ell}$.) This pairing has the following key adjointness property for the Hecke operators,
$$
\la [U\zeta U]f,f'\ra=\la f,[U' I(\zeta)^{-1}U']f'\ra,
$$ 
for any $\zeta \in G(\A_{F^+}^{\infty})$ with $\zeta_{\ell}=1$. As a special case, we see that $T_w^{(j)}$ is self-adjoint -- in the natural sense. More precisely, the adjoint map, which sends $[U\zeta U]$ to $[U' I(\zeta)^{-1}U']$, defines an isomorphism of Hecke algebras,
$$
\text{$\T_{a,\{\rho_v\},\varnothing}^T(U)\overset{\sim}{\longrightarrow} \T_{a,\{\rho_v\},\varnothing}^T(U')'$, $\y$ $T_{w}^{(j)}\mapsto T_{w}^{(j)}$,}
$$
where the latter algebra is defined in the obvious way -- sitting inside the endomorphisms of $S_{a,\{\rho_v\},\varnothing}'(U',A)$. In particular, without further comment, we will view maximal ideals $\m \subset \T_{a,\{\rho_v\},\varnothing}^T(U)$ as maximal ideals of $\T_{a,\{\rho_v\},\varnothing}^T(U')'$, and vice versa. 

\bigskip

Recall from commutative algebra that
$$
\T_{a,\{\rho_v\},\varnothing}^T(U)\overset{\sim}{\longrightarrow} \prod_{\m} \T_{a,\{\rho_v\},\varnothing}^T(U)_{\m},
$$
where $\m$ runs over the finitely many maximal ideals (and similarly for the primed version). Correspondingly, the module $S_{a,\{\rho_v\},\varnothing}(U,\OO)$ decomposes as a direct sum $\oplus_{\m} S_{a,\{\rho_v\},\varnothing}(U,\OO)_{\m}$. Upon tensoring with $\bar{K}$, each $S_{a,\{\rho_v\},\varnothing}(U,\OO)_{\m}$ breaks up into Hecke-eigenspaces, $\oplus _{\p \subset \m} S_{a,\{\rho_v\},\varnothing}(U,\bar{K})_{\p}$. Given these observations, it is easy to see that $\la,\ra$ restricts to a perfect pairing between $\m$-localizations, 
$$
\la,\ra: S_{a,\{\rho_v\},\varnothing}(U,\OO)_{\m} \times S_{a,\{\rho_v\},\varnothing}'(U',\OO)_{\m}\longrightarrow \OO.
$$

\bigskip

Now consider levels $U=U_{\Sigma}U^{\Sigma}$, as in the definition of $H=H_{\Sigma}$. Our immediate goal is to extend the duality to infinite level, after shrinking $U_{\Sigma}\rightarrow \{1\}$. Since $\Sigma \subset S$ consists of banal places, the pro-order $|U_{\Sigma}|$ is invertible in $\OO$, and we may modify the above pairing $\la,\ra$ by suitable constants in order to make them compatible with change of level $U_{\Sigma}$. In the definition of $H$, we had a non-Eisenstein maximal ideal $\m$ occurring at full level $U(\varnothing)$, which we pulled back to each $\T_{a,\{\rho_v\},\varnothing}^T(U_{\Sigma}U^{\Sigma})$. Then,
$$
H=H_{\Sigma}=S_{a,\{\rho_v\},\varnothing}(U^{\Sigma},\OO)_{\m}=\underset{U_{\Sigma}}{\varinjlim}S_{a,\{\rho_v\},\varnothing}(U_{\Sigma}U^{\Sigma},\OO)_{\m}.
$$
Similarly, let
$$
H'=H_{\Sigma}'=S'_{a,\{\rho_v\},\varnothing}(I(U^{\Sigma}),\OO)_{\m}=\underset{U_{\Sigma}}{\varinjlim}S_{a,\{\rho_v\},\varnothing}'(I(U_{\Sigma}U^{\Sigma}),\OO)_{\m},
$$
where we view $\m$ as a maximal ideal of the primed Hecke algebras, via the isomorphism described above. 

The upshot of our discussion is the following duality.

\begin{lem}
There is a (smoothly) perfect pairing $\la,\ra: H \times H' \longrightarrow \OO$, relative to which each $T_w^{(j)}$ is self-adjoint, and such that 
$\la gh,h'\ra=\la h, I(g)^{-1}h'\ra$ for all $g \in G(F_{\Sigma}^+)$. In conclusion, 
$$
\text{$H'\simeq \Hom_{\OO}(H,\OO)^{\infty}$, $\y$ $H\simeq \Hom_{\OO}(H',\OO)^{\infty}$,}
$$
as smooth $G(F_{\Sigma}^+)$-representations defined over $\T_{\Sigma}$. (Intertwining the contragredient action with the action via $I$ on $H'$.)
\end{lem}

The superscript $\infty$ means we take the smooth vectors in $\Hom_{\OO}(-,\OO)$. 

\subsection{A second application of Ihara's lemma}

We already used Ihara's lemma (Conjecture 1) to get $R=\T$ in the non-minimal case, and hence the freeness of $X_{\Sigma}$ over $\T_{\Sigma}$. In this subsection we will apply it -- or rather the stronger version, Conjecture 2 -- in a different way, to see that the generic constituent of $H/\m H$ must in fact be a quotient. 

\bigskip

For a moment, we look at the full $\OO$-linear dual $M=\Hom_{\OO}(H',\OO)$, with its natural $\T_{\Sigma}[G(F_{\Sigma}^+)]$-module structure. With this notation, $H \simeq M^{\infty}$.

\begin{lem}
$M/\m M\simeq \Hom_k((H'/\varpi H')[\m],k)$.
\end{lem}

\noindent {\it{Proof}}. Since $H'$ is $\varpi$-adically complete (where $\lambda=\varpi\OO$) we have
$$
M \overset{\sim}{\longrightarrow} \underset{i}{\varprojlim} \Hom_{\OO}(H'/\varpi^i H', \OO/\varpi^i\OO)=\underset{i}{\varprojlim} M/\varpi^i M.
$$
In particular, for $i=1$, we find that $M/\varpi M\simeq \Hom_k(H'/\varpi H',k)$. Comparing the largest quotients annihilated by $\m$ yields the result. (This roughly follows the proof of Proposition C.5, in Appendix C, on p. 104 of \cite{9} -- and the remarks preceding it.) $\square$

\begin{lem}
Assume Conjecture 2 holds.  Then every irreducible $k[G(F_{\Sigma}^+)]$-quotient of $H/\m H$ is generic.
\end{lem}

\noindent {\it{Proof}}. First off, observe that
$$
(H'/\varpi H')[\m]\simeq S_{a,\{\rho_v\},\varnothing}'(I(U^{\Sigma}),k)_{\m}[\m].
$$
Let $\pi=\otimes_{v\in \Sigma}\pi_v$ be an absolutely irreducible $G'(F_{\Sigma}^+)$-submodule hereof. Then, for every fixed $v \in \Sigma$, there is a $G'(F_v^+)$-equivariant embedding,
$$
\pi_v \hookrightarrow S_{a,\{\rho_v\},\varnothing}'(U'^v,\bar{k})_{\m}[\m],
$$
for some sufficiently small level $U'^v=U_{\Sigma}'^vI(U^{\Sigma})$ away from $v$. However, it may not be the case that $\pi_v$ sits inside a cyclic submodule $\la G'(F_v^+)f\ra$,
for an eigenform $f$ of full level $U_v'=G'(\OO_{F_v^+})$ at $v$ -- in which case Conjecture 1 does not suffice. However, since we appeal to Conjecture 2, we can infer that
$\pi_v$ is generic. We conclude that $\pi_{\Sigma}$ is a generic representation of $G'(F_{\Sigma}^+)$

\bigskip

Finally, taking smooth vectors $(\cdot)^{\infty}$ commutes with extending scalars (by banality of the places in $\Sigma \subset S$), so by the previous lemma,
$$
H/\m H=M^{\infty}/\m M^{\infty}=(M/\m M)^{\infty}\simeq \Hom_k((H'/\varpi H')[\m],k)^{\infty}.
$$
All the smooth representations involved here are admissible. We deduce that 
$$
(H'/\varpi H')[\m] \simeq \Hom_k(H/ \m H,k)^{\infty}.
$$
Consequently, if $\pi_{\Sigma}$ is an irreducible $k[G(F_{\Sigma}^+)]$-quotient of $H/\m H$, its smooth dual 
$\pi_{\Sigma}^{\vee}=\Hom_k(\pi_{\Sigma},k)^{\infty}$ becomes a submodule of $(H'/\varpi H')[\m]$, and is therefore generic as just explained. Taking duals again, 
we see $\pi_{\Sigma}$ itself must be generic. $\square$

\bigskip

This allows us to finish the proof of our main result: 

\begin{thm}
Admit Conjecture 2. Then $H/\m H$ has a unique irreducible quotient. This quotient is absolutely irreducible and generic. All other constituents of $H/\m H$ are non-generic:
In the terminology of \cite{5}, the smooth dual
$$
(H/\m H)^{\vee}\simeq S_{a,\{\rho_v\},\varnothing}'(U'^{\Sigma},k)_{\m}[\m]
$$
is essentially AIG (short for "absolutely irreducible and generic").  
\end{thm}

\noindent {\it{Proof}}. Let $\pi_{\Sigma}$ be an irreducible quotient of $H/\m H$, which is necessarily generic by Lemma 9. Now, by Lemma 6, $H/\m H$ has a unique generic constituent (counting multiplicity) -- which must then be $\pi_{\Sigma}$. As a result, $H/\m H$ has $\pi_{\Sigma}$ as its unique irreducible quotient -- which is generic -- and all other constituents are non-generic. $\square$

\bigskip

For $n=2$ this result is unconditional, and we expect it to play a key role in proving strong local-global compatibility results in the $p$-adic Langlands program for unitary groups in two variables, in the vein of \cite{9}. This is joint work in progress with Chojecki -- in continuation of our "weak" local-global compatibility paper \cite{10}, in which we showed how the $p$-adic local Langlands correspondence for $\GL_2(\Q_p)$ appears in the completed cohomology of the tower of finite sets for $G$.

\subsection{The interpolative property of $H$}

We need to verify that $H$ interpolates the (generic) local Langlands correspondence; or rather its dual $\tilde{\pi}(\cdot)$. (Property (2) in Theorem 6.2.1 in \cite{5}.)
Our multiplicity one assumption ("Hypothesis 1") remains in force.

\begin{lem}
If $\p \subset \T_{\Sigma}$ is a minimal prime, with residue field $\kappa(\p)$, there is a $\kappa(\p)$-linear $G(F_{\Sigma}^+)$-equivariant isomorphism,
$$
\otimes_{v \in \Sigma} \tilde{\pi}(r_{\m,\tilde{v}}\otimes_{\T_{\Sigma}} \kappa(\p)) \overset{\sim}{\longrightarrow} \kappa(\p)\otimes_{\T_{\Sigma}} H.
$$
(Recall that we always identify $G(F_{\Sigma}^+)$ with $\prod_{v\in \Sigma}\GL_n(F_{\tilde{v}})$.)
\end{lem}

\noindent {\it{Proof}}. For simplicity, let us assume $\OO=\T_{\Sigma}/\p$, so that $\kappa(\p)=K$. Then,
$$
\kappa(\p)\otimes_{\T_{\Sigma}} H=\kappa(\p) \otimes_{\T_{\Sigma}/\p}  (H/\p H)= (H\otimes_{\OO}K)/\p (H\otimes_{\OO}K).
$$
Now keep in mind the automorphic description of $H \otimes_{\OO}\bar{K}$ given in the beginning of section 5.1. It breaks up as $\oplus_{\pi}\pi_{\Sigma}$, for a certain family of automorphic representations $\pi$. Moreover, $\T_{\Sigma}$ acts on $\pi_{\Sigma}$ via the eigensystem $\phi_{\pi}: \T_{\Sigma}\rightarrow \bar{K}$. After possibly enlarging $K$,
$$
(H\otimes_{\OO}K)/\p (H\otimes_{\OO}K)=\oplus_{\pi}\pi_{\Sigma}/\p \pi_{\Sigma}.
$$
Now, $\p \pi_{\Sigma}=0$ for $\p=\ker(\phi_{\pi})$, and it is easy to see that $\p \pi_{\Sigma}=\pi_{\Sigma}$ for all other $\p$. Therefore, the above sum reduces to just 
$\pi_{\Sigma}$, where $\pi=\pi(\p)$ now denotes {\it{the}} automorphic representation with $\ker(\phi_{\pi})=\p$. In other words, with
$$
r_{\pi}\simeq r_{\m}\otimes_{\T_{\Sigma}}\kappa(\p).
$$
It remains to note that $\tilde{\pi}(r_{\pi,\tilde{v}})=\pi_v$, at each $v \in \Sigma$, by local-global compatibility away from $\ell$ -- again, under the identification $G(F_v^+)\simeq 
\GL_n(F_{\tilde{v}})$. (See the proof of Proposition 3.3.4 in \cite{1} how to deduce this from the corresponding result for $\GL_n$ -- due to Taylor, Yoshida, Shin, and others.) We should point out why $\pi_v$ is generic. Indeed, $r_{\pi}$ is irreducible (since $\bar{r}_{\pi}\simeq \bar{r}_{\m}$ is; $\m$ being non-Eisenstein), so the base change of $\pi$ to $\GL_n(\A_{F})$ must be cuspidal, and therefore (globally) generic. Thus, in our situation, there is no discrepancy between between local Langlands and "generic" local Langlands. 
$\square$

\subsection{$H^{tf}$ and local Langlands in families}

First recall that an $A$-module $M$ is {\it{torsionfree}} if the map $m \mapsto am$ defines an injection $M \rightarrow M$, for any non-zerodivisor $a \in A$. Equivalently, in a fancier language, every associated prime of $M$ is contained in an associated prime of $A$. When $A$ is reduced, the latter are just the minimal primes. We let $M_{tor}$ denote the subset of $m \in M$ for which $Ann_A(m)$ contains a non-zerodivisor. Note that $M_{tor}$ is an $A$-submodule, because the non-zerodivisors form a multiplicative subset.
We let $M^{tf}=M/M_{tor}$ be the maximal $A$-torsionfree quotient of $M$.

\medskip

\noindent In the next subsection we show that in fact $H$ is $\T_{\Sigma}$-torsionfree. Here we wish to point out that the interpolative property is somehow insensitive to torsion.
Indeed, one can easily deduce the weaker statement that $H^{tf}$ interpolates local Langands (knowing $H$ does).

\begin{lem}
$H \otimes_{\T_{\Sigma}} \kappa(\p)\overset{\sim}{\longrightarrow}H^{tf}\otimes_{\T_{\Sigma}} \kappa(\p)$, for every minimal prime $\p \subset \T_{\Sigma}$.
\end{lem}

\noindent {\it{Proof}}. Let $\II \subset H \otimes_{\T_{\Sigma}}\kappa(\p)$ denote the image of $H_{tor} \otimes_{\T_{\Sigma}} \kappa(\p)$ under the natural map. Clearly 
$G(F_{\Sigma}^+)$ preserves $H_{tor}$, and therefore this $\II$ is a $G(F_{\Sigma}^+)$-invariant subspace. As we have just seen, as part of the proof of the previous Lemma, 
$H \otimes_{\T_{\Sigma}}\kappa(\p)\simeq \pi_{\Sigma}$ is irreducible. Consequently, if $\II \neq 0$, it must be the whole space $H \otimes_{\T_{\Sigma}}\kappa(\p)$. Applying the (composite) functor $\LL_{n,\Sigma}$ from section 5.3,
$$
\LL_{n,\Sigma}(\II)=\LL_{n,\Sigma}(H)\otimes_{\T_{\Sigma}}\kappa(\p)=X_{\Sigma}\otimes_{\T_{\Sigma}}\kappa(\p)=\kappa(\p),
$$
where, in the last step, we have used that $X_{\Sigma}\simeq \T_{\Sigma}$ (the main result from 4.3, which relies on Ihara's lemma). On the other hand,  $\LL_{n,\Sigma}(\II)$
is a quotient of $\LL_{n,\Sigma}(H_{tor})\otimes_{\T_{\Sigma}}\kappa(\p)$, and
$$
\LL_{n,\Sigma}(H_{tor}) \subset \LL_{n,\Sigma}(H)_{tor}=X_{\Sigma,tor}=0,
$$
again using $X_{\Sigma}\simeq \T_{\Sigma}$ in the last step. We deduce that $\LL_{n,\Sigma}(\II)=0$. This is a contradiction. We conclude that $\II=0$. $\square$

\medskip

\noindent {\it{Proof of Theorem 1 -- up to torsion}}: $H^{tf}$ is tautologically $\T_{\Sigma}$-torsionfree, and we have just shown (in the last two lemmas) that it interpolates the local Langlands correspondence. Moreover, there is a $G(F_{\Sigma}^+)$-equivariant embedding,
$$
(H^{tf}/\m H^{tf})^{\vee} \hookrightarrow (H/\m H)^{\vee}.
$$
By our main result, Theorem 4, $(H/\m H)^{\vee}$ is essentially AIG. A fortiori, so is $(H^{tf}/\m H^{tf})^{\vee}$. Theorem 6.2.1 in \cite{5} shows that there is at most one module with these three properties. Having proved its existence, we infer that 
$$
\tilde{\pi}(\{r_{\m,\tilde{v}}\}_{v \in \Sigma})\simeq H^{tf},
$$
where we adapt the notation from 6.2.2 in \cite{5}. It remains to show $H=H^{tf}$.

\subsection{$H$ is torsion-free over $\T_{\Sigma}$}

We finish off by showing that in fact $H$ is $\T_{\Sigma}$-torsionfree.

\begin{lem}
$H=H^{tf}$ (at least after enlarging the coefficient field $K$).
\end{lem}

\noindent {\it{Proof}}. Suppose there is a non-zero element $h \in H_{tor}$. Thus $Ann_{\T_{\Sigma}}(h)$ contains a non-zerodivisor $t \in \T_{\Sigma}$, say. Recall the automorphic description of $H \otimes_{\OO}\bar{K}$ from 5.1,
$$
H \otimes_{\OO}\bar{K} \simeq \bigoplus_{\pi}\pi_{\Sigma},
$$
a finite direct sum (by Lemma 7). As explained at the end of the proof of Proposition 2, each $\pi$ is completely determined by its Hecke eigensystem, which we continue to denote  $\phi_{\pi}: \T_{\Sigma}\rightarrow \bar{K}$. Since $\T_{\Sigma}$ preserves $H$, in fact $\phi_{\pi}$ maps into $\OO_{\bar{K}}$. We enlarge $K$, if necessary, to arrange for $\OO=\OO_K$ to accommodate the images of all the (finitely many) $\phi_{\pi}$. Via the isomorphism,
$$
\T_{\Sigma}\otimes_{\OO} \bar{K} \simeq \bar{K}^{\Hom_{\OO-alg.}(\T_{\Sigma},\bar{K})}=\bar{K}^{\Hom_{\OO-alg.}(\T_{\Sigma},\OO)},
$$
$t\otimes 1$ corresponds to the tuple $(\phi_{\pi}(t))_{\pi}$. Moreover, $h\otimes 1$ corresponds to some tuple, which we will label $(h_{\pi})_{\pi}$, for some vectors
$h_{\pi}\in \pi_{\Sigma}$. Now, the assumption that $th=0$ translates into $\phi_{\pi}(t)h_{\pi}=0$ -- for all $\pi$. Since $h$ is nonzero, $h_{\pi}\neq 0$ for at least one $\pi$, which we fix for the rest of this proof. For this specific $\pi$, we deduce that $\phi_{\pi}(t)=0$. Now let $e_{\pi}\in \T_{\Sigma}\otimes_{\OO} \bar{K}$ be the element which corresponds to the "standard" vector $(0,\ldots,1,\ldots,0)$, which is $1$ in the $\pi$-slot, and $0$ elsewhere. Then obviously $e_{\pi}(t\otimes 1)=0$ -- which shows $t\otimes 1$ {\it{is}} a 
zero-divisor in the bigger ring $\T_{\Sigma}\otimes_{\OO} \bar{K}$. However, it may not be the case that $e_{\pi} \in \T_{\Sigma}$. Nevertheless, we claim that 
$ae_{\pi}\in \T_{\Sigma}$ for some nonzero $a \in \OO$. To see this, observe that the natural inclusion
$$
i: \T_{\Sigma}\hookrightarrow \OO^{\Hom_{\OO-alg.}(\T_{\Sigma},\OO)}=\OO \times \cdots \times \OO
$$
becomes an isomorphism after $(-)\otimes_{\OO}\bar{K}$. Thus $\cok(i)\otimes_{\OO}\bar{K}=0$, which is to say $\cok(i)$ is $\OO$-torsion (since $\bar{K}$ is flat over $\OO$). 
$e_{\pi}$ does lie in the target of $i$, and so some nonzero multiple $ae_{\pi}$ lies in its image $\im(i)$. In conclusion, we get an element in $\T_{\Sigma}$, which we will continue to denote $ae_{\pi}$, such that $ae_{\pi}t=0$. This shows $t$ is a zerodivisor in $\T_{\Sigma}$, contrary to our assumption. $\square$

\section{$H/\m H$ and local Langlands mod $\ell$}

We end this paper with a few observations on how $H/\m H$ relates to the mod $\ell$ local Langlands correspondence of Emerton-Helm (as defined in Chapter 5 of \cite{5}).
By our main result (Theorem 1 in the introduction) we know that Ihara's lemma is {\it{equivalent}} to the isomorphism $\tilde{\pi}(\{r_{\m,\tilde{v}}\}_{v \in \Sigma})\simeq H$. 
Reducing mod $\m$,
$$
\tilde{\pi}(\{r_{\m,\tilde{v}}\}_{v \in \Sigma})\otimes_{\T_{\Sigma}} k\simeq H/\m H.
$$
Now, Conjecture 6.2.9 on p. 52 in \cite{5} predicts that there should be an equivariant {\it{surjection}},
$$
\otimes_{v\in \Sigma} \tilde{\bar{\pi}}(\bar{r}_{\m,\tilde{v}})\twoheadrightarrow \tilde{\pi}(\{r_{\m,\tilde{v}}\}_{v \in \Sigma})\otimes_{\T_{\Sigma}} k.
$$
We prove this below -- contingent on Ihara's lemma. Conversely, dualizing such a surjection, $(H/\m H)^{\vee}$ sits as a subrepresentation of 
$\otimes_{v\in \Sigma} \bar{\pi}(\bar{r}_{\m,\tilde{v}})$, which is essentially AIG. Therefore $(H/\m H)^{\vee}$ is itself essentially AIG. Identifying it with a space of mod $\ell$ modular forms, using Petersson duality as in Theorem 4 above, one can {\it{deduce}} Ihara's lemma. Thus we end up with following mod $\ell$ reformulation of Ihara:

\begin{cor}
Ihara's lemma (Conjecture 1 in the introduction) for $G$ is equivalent to the existence of an equivariant surjection,
$$
\otimes_{v\in \Sigma} \tilde{\bar{\pi}}(\bar{r}_{\m,\tilde{v}})\twoheadrightarrow H/\m H.
$$
(Dually, $S'_{a,\{\rho_v\},\varnothing}(U'^{\Sigma},k)_{\m}[\m]\hookrightarrow \otimes_{v\in \Sigma} \bar{\pi}(\bar{r}_{\m,\tilde{v}})$, with notation as in Theorem 4.
In particular,
$$
\otimes_{v\in \Sigma} \bar{\pi}(\bar{r}_{\m,\tilde{v}}) \hookrightarrow env (H/\m H)^{\vee},
$$
where $env$ denotes the essentially AIG envelope -- introduced in Definition 3.2.6 on p. 21 in \cite{5}.)
\end{cor}

\noindent {\it{Proof}}. Let us first assume $\Sigma=\{v\}$ is a singleton, and at the end reduce the general case to this special case. Our goal is to produce a surjection,
$$
\tilde{\bar{\pi}}(\bar{r}_{\m,\tilde{v}})\twoheadrightarrow \tilde{\pi}(r_{\m,\tilde{v}})\otimes_{\T_{\Sigma}}k.
$$
This would follow immediately if we knew $\bar{\pi}(\cdot)$ has the strong property (2') in Remark 1.5.2 of \cite{5} -- but we don't. To make use of the weaker property (2) in Theorem 1.5.1 of loc. cit., we choose some minimal prime $\p \subset \T_{\Sigma}$ and look at the deformation of $\bar{r}_{\m,\tilde{v}}$ given by specializing $r_{\m,\tilde{v}}$ at $\p$,
$$
r_{\m,\tilde{v}}\otimes_{\T_{\Sigma}} \T_{\Sigma}/{\p}: \Gamma_{F_{\tilde{v}}} \longrightarrow \GL_n(\T_{\Sigma}/{\p}).
$$
This is a representation over $\T_{\Sigma}/\p$, which we may identify with just $\OO$ after possibly enlarging our coefficient field $K/\Q_{\ell}$ -- as we did in the proof of Lemma 14 (except that we don't pass to the residue field $\kappa(\p)$). By property (2), alluded to earlier, we now {\it{know}} there is a surjection,
$$
\tilde{\bar{\pi}}(\bar{r}_{\m,\tilde{v}})\twoheadrightarrow \tilde{\pi}(r_{\m,\tilde{v}}\otimes_{\T_{\Sigma}} \T_{\Sigma}/{\p})\otimes_{\T_{\Sigma}/{\p}}k.
$$
Moreover, by Lemma 6.3.12 on p. 60 in \cite{5}, we know that 
$$
\tilde{\pi}(r_{\m,\tilde{v}}\otimes_{\T_{\Sigma}} \T_{\Sigma}/{\p})\simeq (\tilde{\pi}(r_{\m,\tilde{v}})\otimes_{\T_{\Sigma}} \T_{\Sigma}/{\p})^{\T_{\Sigma}/{\p}-tf},
$$
the maximal $\T_{\Sigma}/\p$-torsionfree quotient. Once we show $\tilde{\pi}(r_{\m,\tilde{v}})\otimes_{\T_{\Sigma}} \T_{\Sigma}/{\p}$ is already torsionfree, we are done.
Indeed, we would then get a map from $\tilde{\bar{\pi}}(\bar{r}_{\m,\tilde{v}})$ onto
$$
\tilde{\pi}(r_{\m,\tilde{v}}\otimes_{\T_{\Sigma}} \T_{\Sigma}/{\p})\otimes_{\T_{\Sigma}/{\p}}k\simeq \tilde{\pi}(r_{\m,\tilde{v}})\otimes_{\T_{\Sigma}} \T_{\Sigma}/{\p}\otimes_{\T_{\Sigma}/{\p}}k\simeq \tilde{\pi}(r_{\m,\tilde{v}}) \otimes_{\T_{\Sigma}}k,
$$
as desired. To show $H/\p H$ is $\T_{\Sigma}/\p$-torsionfree, we tweak the proof of Lemma 16 a bit. Recall from the proof of Lemma 14 that
$$
H/\p H \otimes_{\OO}\bar{K}=(H \otimes_{\OO}\bar{K})/\p(H \otimes_{\OO}\bar{K})=\pi_{\Sigma},
$$
where $\pi=\pi(\p)$ is the automorphic representation whose eigensystem $\phi_{\pi}$ has kernel $\p$. Thus $t \in \T_{\Sigma}$ acts on $H/\p H$ via multiplication by $\phi_{\pi}(t)$. If $t$ annihilates some nonzero element of $H/\p H$, we must have $\phi_{\pi}(t)=0$ -- which shows $H/\p H$ is torsionfree, and concludes the proof of the Corollary in the singleton-case.

\medskip

It remains to deal with the general case, where $\Sigma$ is not necessarily a singleton. Taking the tensor product $\otimes_{v \in \Sigma}$ (over $k$) of all the surjections defined in the previous paragraph, we get a map
$$
\otimes_{v\in \Sigma} \tilde{\bar{\pi}}(\bar{r}_{\m,\tilde{v}})\twoheadrightarrow \otimes_{v\in \Sigma}(\tilde{\pi}(r_{\m,\tilde{v}})\otimes_{\T_{\Sigma}}k)\simeq
(\otimes_{v \in \Sigma}\tilde{\pi}(r_{\m,\tilde{v}}))\otimes_{\T_{\Sigma}} k,
$$
(where the tensor product $\otimes_{v \in \Sigma}\tilde{\pi}(r_{\m,\tilde{v}})$ obviously is over $\T_{\Sigma}$). Compose it with 
$$
\otimes_{v \in \Sigma}\tilde{\pi}(r_{\m,\tilde{v}}) \twoheadrightarrow (\otimes_{v \in \Sigma}\tilde{\pi}(r_{\m,\tilde{v}}))^{tf}\simeq \tilde{\pi}(\{r_{\m,\tilde{v}}\}_{v \in \Sigma})
$$
after tensoring with $(-)\otimes_{\T_{\Sigma}}k$. The last isomorphism above is the content of Proposition 6.2.4 on p. 51 in \cite{5}. This proves the Corollary.
$\square$



\noindent {\sc{Department of Mathematics, UCSD, La Jolla, CA, USA.}}

\noindent {\it{E-mail address}}: {\texttt{csorensen@ucsd.edu}}

\end{document}